\documentclass[a4paper,draft]{amsproc}
\usepackage{amssymb}
\usepackage[hyphens]{url} 
\theoremstyle{plain}
 \newtheorem{thm}{Theorem}[section]
 
 \newtheorem{lem}{Lemma}[section]
 \newtheorem{conj}{Conjecture}[section]
 
\theoremstyle{definition}
 \newtheorem{exm}{Example}[section]
 
\theoremstyle{remark}
 \newtheorem{rem}{Remark}[section]
 \numberwithin{equation}{section}

\renewcommand{\leq}{\leqslant}
\renewcommand{\ge}{\geqslant}\renewcommand{\geq}{\geqslant}

\begin{document}
%------------------------------------------------------
%------------------------------------------------------
%--- Mettre les exercices ici---
%--- Les corrections sont automatiquement ajoutees ---
%%%%%%%%%%%%%%   exercice &1    %%%%%%%%%%%%%%%%%%%%%%%%%%%%%%%%%%%%%%%%%%%%%%%%%%%%%%%%%%%%
\title[Diophantine approximation and continued fraction expansion]{Diophantine approximation and continued fraction expansion for quartic power series over $\mathbb{F}_{3}$}
%\subtitle{Diophantine approximation}{if too long for running head}

\author{Khalil Ayadi, Awatef Azaza and Salah Beldi}

\maketitle

\textbf{Abstract} While Roth's theorem states that the irrationality measure of all the irrational algebraic numbers is 2, and the same holds true
over function fields in characteristic zero, some counter-examples were found over function fields in positive characteristic.
This was put forward first by Mahler in 1949, in his fundamental paper on Diophantine approximation \cite{M}. It seems
that, except for particular elements, as power series with bounded partial quotients, Roth's theorem holds. Until now, only one element,
with unbounded partial quotients, discovered by  Mills and  Robbins \cite{MR} in 1986,
has been recognized having this property.
It concerns a quartic power series over $\mathbb{F}_{3}$ having a continued fraction expansion with remarkable pattern. This continued
fraction expansion was explicitly described by  Buck and Robbins \cite{BR}, and later by Lasjaunias \cite{LA2} who used
another method somewhat easier. Furthermore, Lasjaunias \cite{LA2} improve the value of its irrationality measure in relation with Roth's theorem.
We will see that this  power series is included in a large quartic power series family, for which the continued
fraction expansion and the irrationality measure can be explicitly given. Moreover, we will study the rational approximation of other
examples of quartic power series  over $\mathbb{F}_{3}$ and we will extend the set of counter-examples initiated by Mahler.\\
\textbf{keywords:}Finite fields, Formal power series, Continued fraction\\
\textbf{MSC:}11J61, 11J70

\maketitle

\section{Introduction}
Let $p$ be a prime number and let $\mathbb{F}$ be a finite field of characteristic $p$.
We let $\mathbb{F}[T]$, $\mathbb{F}(T)$ and $\mathbb{F}((T^{-1}))$ respectively denote, the ring of polynomials, the field of rational
functions and the field of power series in $1/T$ over $\mathbb{F}$, where $T$ is a formal indeterminate. These
fields are valuated by the ultrametric absolute value introduced on $\mathbb{F}(T)$ by
$|P/Q|= |T|^{\deg(P)-\deg(Q)},$ where $|T|>1$ is a fixed real number. We recall that each irrational
(rational) element $\alpha$ of $\mathbb{F}((T^{-1}))$ can be expanded as an infinite (finite) continued fraction. This will be
denoted $\alpha=[a_{0}, a_{1}, . . . , a_{n}, . . . ]$ where the $a_{i} \in \mathbb{F}[T]$, with $\deg(a_{i})>0$ for $i\geq1$, are the partial
quotients and the tail $\alpha_{i} = [a_{i}, a_{i+1}, . . . ] \in \mathbb{F}((T^{-1}))$ is the complete quotient. As in
the classical theory, we define recursively the two sequences of polynomials $(P_{n})_{n\geq0}$
and $(Q_{n})_{n\geq0}$ by
$P_{n}= a_{n}P_{n-1} + P_{n-2}$ and $Q_{n}= a_{n}Q_{n-1} + Q_{n-2}$,
with the initial conditions $P_{0}= a_{0}$, $P_{1}= a_{1}a_{2} + 1$, $Q_{0}= 1$ and $Q_{1}= a_{2}$. We
have $P_{n+1}Q_{n}- Q_{n+1}P_{n}= (-1)^{n}$, whence $P_{n}$ and $Q_{n}$ are coprime polynomials. The
rational function $P_{n}/Q_{n}$ is called a convergent to $\alpha$ and we have $P_{n}/Q_{n}= [a_{0}, a_{1}, ..., a_{n}]$ and
$P_{n}/P_{n-1}= [a_{n}, a_{n-1}, ..., a_{0}]$.
It is easily checked the following property of continued fraction: when $B, C$ are nonzero polynomials in $\mathbb{F}[T]$, then
\begin{eqnarray}\label{sd}C[Ba_{0},Ca_{1},Ba_{2}, . . .] = B[Ca_{0},Ba_{1},Ca_{2}, . . .].\end{eqnarray}
As for real numbers, the continued fraction expansion of formal power series is fundamental to measure the quality of their rational approximation.
The irrationality measure(or the approximation exponent) of an irrational power series $\alpha\in \mathbb{F}((T^{-1}))$ is defined by:
$$\nu(\alpha)=-\displaystyle\limsup_{|Q|\longrightarrow\infty}\log(|\alpha -P/Q|)/\log(|Q|)$$
where $P, Q\in \mathbb{F}[T]$. It is directly related to the growth of
the sequence of the degrees of the partial quotients in the continued fraction
expansion of $\alpha$. Indeed we have
\begin{eqnarray}\label{A}\nu(\alpha)= 2 + \displaystyle\limsup_{n>1}(\deg(a_{n+1})/\sum_{1\leq i\leq n}\deg(a_{i})).\end{eqnarray}
Note that the irrationality measure is stable under a M\"{o}bius transformation of nonzero determinant.\\
For a general presentation of continued fractions and diophantine approximation in the function
field case, the reader may consult \cite{S} or (\cite{dt} Chap. 9).\\
We have to be concerned with infinite continued fractions in $\mathbb{F}((T^{-1}))$ which are algebraic over $\mathbb{F}(T)$.
The study of their rational approximation was initiated by  Mahler \cite{M}. The starting point in the study of rational approximation to algebraic real
numbers is a famous theorem of Liouville established in 1850. This theorem was adapted by Mahler in the fields of power series with an arbitrary base field: if $\alpha$ is an element of $\mathbb{F}((T^{-1}))$,
algebraic of degree $n>1$ over $\mathbb{F}(T)$, then for all element $P/Q$ of $\mathbb{F}(T)$,  there exists a positive real number $c$ such that
$$|\alpha- P/Q|\geq c/|Q|^{n}.$$
This result implies that $\nu(\alpha)\leq n$. In the case of real numbers, it is well known that Liouville's theorem has been
improved until Roth's theorem was established. In 1955, Roth \cite{R} proved that for any irrational algebraic real number $\alpha$, $\nu(\alpha)=2$.
These improvements on the exponents have many applications to solve Diophantine equations and transcendental questions. It is transposed
in fields of power series if the base field has the characteristic zero
as proved by Uchiyama in 1960 \cite{U}, the exponents of irrational algebraic power series is still 2. In this case the exponent $n$ in the right
hand side of the above inequality can be replaced by $2+\epsilon$ for all $\epsilon>0$
. But a naive analog of Roth's theorem now fails in positive characteristic and consequently the
study of rational approximation to algebraic elements becomes more complex. Mahler \cite{M} gave an example showing that the approximation exponent
$\nu(\alpha)$ could reach  the $\alpha$'s degree $n$. He has considered the irrational solution in $\mathbb{F}_{p}((T^{-1}))$  of
the equation $x = 1/T + x^{p}$. For this element $\alpha$, algebraic of degree $p$,
we have rationals $P/Q$, with $|Q|$ arbitrarily large, and $|\alpha - P/Q| = |Q|^{-p}$.\\
Regarding diophantine approximation and continued fractions, a particular subset of elements in $\mathbb{F}((T^{-1}))$,
algebraic over $\mathbb{F}(T)$ is worth considering. Let $r=p^{t}$
with $t\geq 0$, we denote by
$\mathcal{H}(r)$ the subset of irrational $\alpha$ belonging to $\mathbb{F}((T^{-1}))$ and satisfying an algebraic equation of the particular
form $A\alpha^{r+1} + B\alpha^{r}
+ C\alpha + D = 0$, where $A, B, C$ and $D$ belong to $\mathbb{F}[T]$. Note that $\mathcal{H}(1)$
is simply the set of quadratic irrational elements in $\mathbb{F}((T^{-1}))$. The union of the subsets $\mathcal{H}(p^{t})$, for $t \geq0$, denoted by
$\mathcal{H}$, is the set of hyperquadratic power series.\\
The rational approximation properties of the elements of $\mathcal{H}$, were
studied independently by  Voloch \cite{V}, and  de Mathan \cite{dM}. They proved that:\\
\emph{If $\alpha\in \mathcal{H}$, and $P/Q \in \mathbb{F}(T)$, either we have}
\begin{eqnarray}\label{b1}\displaystyle\liminf_{|Q|\longrightarrow\infty}|Q|^{2}|\alpha -P/Q|>0 \end{eqnarray}
\emph{or there exists a real number $\mu>2$ such that }
\begin{eqnarray}\label{b2}\displaystyle\liminf_{|Q|\longrightarrow\infty}|Q|^{\mu}|\alpha -P/Q|<\infty.\end{eqnarray}
 With
respect to this, de Mathan and Lasjaunias \cite{LA3}, have shown that if an
algebraic element does not belong to $\mathcal{H}$, then it cannot be too well aproximated by
rationals : if $\alpha \not\in \mathcal{H}$ and it is algebraic of degree $n > 1$ over $\mathbb{F}(T)$, then, for all $\epsilon> 0$,
we have $|\alpha- P/Q|> |Q|^{-([n/2]+1+\epsilon)}$
, for all $P/Q \in \mathbb{F}(T)$ with $|Q|$ large enough. This last property highlights the peculiarity of the set $\mathcal{H}$.
If rational approximation to certain hyperquadratic power series is
well known, this is also due to the possibility of describing explicitly their
continued fraction expansion. The first works in this area were undertaken
by Baum and Sweet \cite{BS}. Later this has been done for many examples and
for different subclasses of hyperquadratic elements (see in particular
\cite{S}). Nevertheless, the possibility of describing the continued fraction expansion for all hyperquadratic
power series is yet an open problem. In \cite{MR} Mills and Robbins
studied this problem by describing an algorithm to obtain, in
certain cases, the continued fraction expansion for an hyperquadratic power
series. They ultimately considered $(p. 403)$ the following
algebraic equation: \begin{eqnarray}\label{1}x^{4}+ x^{2}- Tx+ 1 = 0.\end{eqnarray}
They observed that this equation
has a unique solution $\alpha$ in $\mathbb{F}_{p}((T^{-1}))$ for all primes $p$ noting that for this
solution, the continued fraction expansion has a remarkable pattern in both
cases $p= 3$ and $p= 13$. The expansion in the case $p= 3$ was explicitly described by  Buck and  Robbins \cite{BR}. Indeed, they
recursively defined the following polynomial sequences:
$$\Omega_{0}=\emptyset,\ \ \ \Omega_{1}=T, \ \ \ \Omega_{n}=\Omega_{n-1}, -T, \Omega^{(3)}_{n-2}, -T,\Omega_{n-1}\ \ \ for\ \ n\geq2.$$
(here $\Omega^{(3)}_{k}$ denotes the sequence obtained by cubing each element of $\Omega_{k}$ and commas indicating juxtaposition of sequences); then they
proved that $[0, \Omega_{n}]$ is the beginning for all $n>0$ of the continued fraction expansion of this solution. This element satisfies, $\displaystyle\liminf_{|Q|\longrightarrow\infty}|Q|^{2}|\alpha -P/Q|=0$
and $\displaystyle\liminf_{|Q|\longrightarrow\infty}|Q|^{\mu}|\alpha -P/Q|=\infty$ for all $\mu>2$.
So it satisfies neither \eqref{b1} nor \eqref{b2}. Thus it does not belong to the set $\mathcal{H}$.
This result was given by Lasjaunias in \cite{LA2}, by proving that there are two real positive constants $\lambda_{1}$ and
$\lambda_{2}$ such that, for some rationals $P/Q$ with $|Q|$ arbitrary large, we have $|\alpha -P/Q|\leq|Q|^{-(2 +\lambda_{1}\sqrt{\log|Q|})}$, and for
all rationals $P/Q$ with $|Q|>1$, we have $|\alpha -P/Q|\geq|Q|^{-(2 +\lambda_{2}\sqrt{\log|Q|})}$. For instance, this element seems to be the first
algebraic element for which the exponent approximation is equal to 2, although its partial quotients are unbounded.\\
Note that for each prime $p> 3$, the continued fraction expansion of the solution of \eqref{1} is remarkable and
it has two different regular patterns and two different values of irrationality measure according to the remainder, $1$ or $2$, in the division of $p$ by $3$,
see \cite{L5} and \cite{L2} for more details.\\

Our work is organized as follow. In the second Section we will compute the continued fraction and the approximation exponent of some quartic power series
which are hyperquadratic over $\mathbb{F}_{3}$. For this, we will use an earlier Theorem which allows us to
determine the approximation exponent of algebraic element when it is large enough, i.e, not close to $2$. The basic idea of this Theorem is due to
 Voloch \cite{V}. It has been improved by de Mathan \cite{dM1}.
\begin{thm}\label{t}$($\cite{L1} p. $219)$ Let $\alpha\in \mathbb{F}((T^{-1}))$. Assume that there is a sequence\\
$(P_{n}, Q_{n})_{n\geq0}$, with $P_{n}$, $Q_{n}\in \mathbb{F}[T]$, satisfying the following conditions:\\
$(1)$ There are two real constants $\lambda>0$ et $\mu>1$, such that
$$|Q_{n}|=\lambda|Q_{n-1}|^{\mu} \ \ \ and\ \ \ |Q_{n}|>|Q_{n-1}|\ \ for\ \  all\ \  n\geq1.$$
$(2)$ There are two real constants $\rho>0$ and $\gamma>1+\sqrt{\mu}$, such that
$$\bigg|\alpha - \displaystyle\frac{P_{n}}{Q_{n}}\bigg|=\rho|Q_{n}|^{-\gamma}\ \ \ for\ \  all\ \ n\geq0.$$
Then we have $\nu(\alpha)=\gamma$.
\end{thm}
This Theorem allows us to find the approximation exponent of several examples of hyperquadratic elements(see \cite{kh}, \cite{L1}).\\
In Section $3$ of this work, we will study the continued fraction expansion  of the solution $\alpha$  of the quartic equation
$$C^{2}\alpha^{4}+2C\alpha^{2}-A^{2}\alpha +1=0\ \ \ (E_{1})$$
where $A$ and $C$ are nonzero polynomials in $\mathbb{F}_{3}[T]$ such that $A$ is not constant, $C$ divides $A$ and $\deg A\geq\deg C$.
By computing the approximation exponent of the solution of this equation, we will prove that is not-hyperquadratic.
Our observation, based on computer calculation giving a finite number of partial quotients for many couples $(A, C)$ of polynomials, implies that
the solution of the equation $(E_{1})$, has very regular pattern in its continued fraction expansion.
Note that this equation can be viewed as a generalization
of the equation \eqref{1} introduced by Mills and Robbins. The properties of rational approximation of $\alpha$ were studied by Lasjaunias in
\cite{LA2} for the case
$A=T$ and $C=-1$. For this case, the tools used to obtain a proof might be well be applied in the general case,
but we are aware that a different approach would be desirable. We will recall the steps of the proof and we will just give our result conjecturally.
Thus we expose a large family of algebraic power series having an approximation exponent value
equal to $2$, even though the degrees of their partial quotients are unbounded.
The great interest of our equation will be to give us the opportunity to introduce and
to describe this family.
\section{Diophantine approximation for some hyperquadratic power series of degree four over $\mathbb{F}_{3}(T)$}
In this section we will study respectively the properties of rational approximations of the solutions of the equations
$$C\beta^{4}-A\beta+1=0 \ \ \ \ (W_{1})\ \  and\ \  -\beta^{4}-A\beta+C=0 \ \ \ \ (W_{2})$$
where $\deg A\geq \deg C$ for $(W_{1})$ and $\deg A> \deg C$ for $(W_{2})$.
\begin{thm}\label{t3} Let $\beta$ be the irrational solution of the equation $(W_{1})$ such that $|\beta|<1$. Assume that $C$ divides $A$.
Then the continued fraction expansion of $\beta$ is $$[b_{0}, b_{1}, \ldots, b_{n},\ldots]$$
such that $b_{0}=0$,  $b_{1}=A$ and for all $n\geq2$:

\begin{eqnarray}\label{ZQ}b_{n}=\left\{%
\begin{array}{ll}
    -Cb_{n-1}^{3} & \hbox{if $n$ is odd ;} \\
    b_{n-1}^{3}/-C  & \hbox{if $n$ is even.} \\
\end{array}%
\right.    \end{eqnarray}
Furthermore, $\nu(\beta)=4$.
\end{thm}
\begin{proofname}. We have $|\beta|<1$ then $b_{0}=0$. Let $\beta_{1}=\beta^{-1}$ then $\beta_{1}$ satisfies the equation
$\beta_{1}^{4}-A\beta_{1}^{3}+C=0$.
Clearly $[\beta_{1}]=b_{1}=A$. In fact, as $|\beta_{1}|>1$ then
$|\beta_{1}^{4}|=|A\beta_{1}^{3}+C|=|A\beta_{1}^{3}|$ so $|\beta_{1}|=|A|$, and since $|\beta_{1}-A|=|C/\beta_{1}^{3}|<1$ then we obtain that
$[\beta_{1}]=A$.
We can write the equation satisfied by $\beta_{1}$ as $\beta_{1}^{3}=\displaystyle\frac{-C}{\beta_{1} -A}$.
So
\begin{eqnarray}\label{sc}\beta_{1}^{3}=-C\beta_{2}.\end{eqnarray}
Applying the Frobenius automorphism to both terms of the identity $\beta_{1}=b_{1}+ 1/\beta_{2}$ and using $\beta_{2}=b_{2}+ 1/\beta_{3}$ we obtain
$b_{2}+\displaystyle\frac{1}{\beta_{3}}=\frac{b_{1}^{3}}{-C} +\frac{1}{-C\beta_{2}^{3}}$. As $C$ divides $A=b_{1}$ then $C$ divides $b_{1}^{3}$, so we get that $b_{2}=b_{1}^{3}/-C$ and
$$\beta_{3}=-C\beta_{2}^{3}.$$ Again, this gives that $b_{3}+\displaystyle\frac{1}{\beta_{4}}= -Cb_{2}^{3}+\frac{-C}{\beta_{3}^{3}}$. So we obtain
$b_{3}=-Cb_{2}^{3}$ and \begin{eqnarray}\label{sck}\beta_{4}=\displaystyle\frac{\beta_{3}^{3}}{-C}.\end{eqnarray}
This gives that $C$ divides $b_{3}$ and \eqref{sc} has the same shape as \eqref{sck}. We now claim that for all $k\ge 1$,
\begin{eqnarray}\label{sckcs}
\left\{%
\begin{array}{ll}
    b_{2k}=b_{2k-1}^{3}/-C, b_{2k+1}=-Cb_{2k}^{3} \\
    \beta_{2k+2}=\beta_{2k+1}^3/-C,\beta_{2k+1}=-C\beta_{2k}^3 \\
\end{array}%
\right.
\end{eqnarray}
Clearly \eqref{sckcs} is true for $k=1$. So we assume \eqref{sckcs} for $k=l\ge 1$. Then
$$\beta_{2l+2}=((b_{2l+1}^3/-C)+\frac{1}{-C\beta_{2l+2}^3}).$$
From \eqref{sckcs} we have $C$ divides $b_{2l+1}^3$. This implies that $b_{2l+2}=b_{2l+1}^3/-C$ and $\beta_{2l+3}=-C\beta_{2l+2}^3$. Then
$$\beta_{2l+3}=-C(b_{2l+2}^3+\frac{1}{\beta_{2l+3}^3})=-Cb_{2l+2}^3+\frac{-C}{\beta_{2l+3}^3},$$
which implies $b_{2l+3}=-Cb_{2l+2}^3$ and $\beta_{2l+4}=\beta_{2l+3}^3/-C$. Thus \eqref{sckcs} is also true for $k=l+1$.
By induction, we see that \eqref{sckcs} holds for all $k\ge 1$.\\
Furthermore, we can verify that the equality \eqref{ZQ} gives that for all $n\geq1$: $$b_{n}=(-1)^{n-1}A^{3^{n-1}}C^{-\frac{3^{n-1}+(-1)^{n}}{4}}.$$
Thus the continued fraction expansion of $\beta$ can be written as
$$[0,A, -A^{3}C^{-1},A^{3^{2}}C^{-2}, -A^{3^{3}}C^{-7},\ldots, (-1)^{n-1}A^{3^{n-1}}C^{-\frac{3^{n-1}+(-1)^{n}}{4}}, \ldots].$$
Now let $a=\deg A$ and $c=\deg C$. Knowing all the partial quotients of $\beta$, we can compute its approximation exponent by the formula \eqref{A}:
\begin{eqnarray*}\nu(\beta)&=&2+\displaystyle\limsup \frac{3^{n}a
-\frac{3^{n}+(-1)^{n+1}}{4}c}{\sum_{k=1}^{n}(3^{k-1}a-\frac{3^{k-1}+(-1)^{k}}{4}c)}\\&=&2 +2=4.\end{eqnarray*}
\end{proofname}
In the next Theorem, we will give the value of $\nu(\beta)$ for $\beta$ satisfying the equation $(W_{1})$ with the condition on the coefficients of
this equation that is: $C$ does not divides $A$.
\begin{thm}\label{t4} Let $\beta$ be the irrational solution of equation $(W_{1})$ such that $|\beta|<1$. Assume that $C$ does not divide $A$.
Then
 $$\nu(\beta)=4-\displaystyle\frac{\deg C}{\deg A}.$$
\end{thm}
\begin{proofname}.
Let $\beta_{1}$ and $\beta_{2}$ be the first and the second complete quotient of $\beta$. So $\beta_{1}$  satisfies
 the equation $\beta_{1}^{4}-A\beta_{1}^{3}+C=0$. We have that $[\beta_{1}]=A$ and since $\beta_{1}=A+ 1/\beta_{2}$ then we can easily see that
 $\beta_{2}$ satisfies the equation $C\beta_{2}^{4}+A^{3}\beta_{2}^{3}+1=0.$
Hence $|\beta_{1}|=|A|$ and $|\beta_{2}|=|A^{3}/C|$. Let $s$ be a positive rational number such that $|A|=|C|^{s}$.
We consider the following sequence: $P_{0}=1$, $Q_{0}=A$ and for $n\geq1$$$P_{n}=Q_{n-1}^{3}$$$$Q_{n}= AQ_{n-1}^{3}-CP_{n-1}^{3}.$$
Then for all $n\geq0$:\\
$$\bigg|\beta-\displaystyle\frac{P_{n}}{Q_{n}}\bigg|=\bigg|\frac{1}{C\beta^{3}-A}-\frac{Q_{n-1}^{3}}{AQ_{n-1}^{3}-CP_{n-1}^{3}}\bigg|=
\bigg|\frac{CP_{n-1}^{3}-CQ_{n-1}^{3}\beta^{3}}{(C\beta^{3}-A)(AQ_{n-1}^{3}-CP_{n-1}^{3})}\bigg|.$$
As $|C\beta^{3}-A|=|A|$ and $|AQ_{n-1}^{3}-CP_{n-1}^{3}|=|AQ_{n-1}^{3}|$ for all $n\geq1$, then we get
$$\bigg|\beta-\displaystyle\frac{P_{n}}{Q_{n}}\bigg|=\frac{|C||P_{n-1}-Q_{n-1}\beta|^{3}}{|A|^{2}|Q_{n-1}|^{3}}=
\frac{|C|}{|A|^{2}}\bigg|\beta-\displaystyle\frac{P_{n-1}}{Q_{n-1}}\bigg|^{3}.$$
 We show by recursion that for all $n\geq0$:
$$\bigg|\beta-\displaystyle\frac{P_{n}}{Q_{n}}\bigg|=\frac{|C|^{\frac{(3^{n}-1)}{2}}}{|A|^{3^{n}-1}}\bigg|\beta-\frac{P_{0}}{Q_{0}}\bigg|^{3^{n}}.$$
\noindent Since $\bigg|\beta-\displaystyle\frac{P_{0}}{Q_{0}}\bigg|=\bigg|\beta-\frac{1}{A}\bigg|=
\frac{|A-\beta_{1}|}{|\beta_{1}||A|}=\frac{1}{|A|^{2}|\beta_{2}|}=\frac{|C|}{|A|^{5}}$ then
$\bigg|\beta-\displaystyle\frac{P_{0}}{Q_{0}}\bigg|^{3^{n}}=|C|^{3^{n}}|A|^{-5.3^{n}}$. So
$$\bigg|\beta-\displaystyle\frac{P_{n}}{Q_{n}}\bigg|=|C|^{\frac{3^{n}-1}{2}}|C|^{3^{n}}|A|^{-(3^{n}-1)}|A|^{-5.3^{n}}=|C|^{\frac{3.3^{n}-1}{2}}
|A|^{-(6.3^{n}-1)}.$$
Let $|A|=|C|^{s}$. Then $\bigg|\beta-\displaystyle\frac{P_{n}}{Q_{n}}\bigg|=|C|^{-\frac{(4s-1)3^{n+1}-2s+1}{2}}.$
Secondly, we have for all  $n\geq1$ $Q_{n}= AQ_{n-1}^{3}-CP_{n-1}^{3}$ then $$|Q_{n}|=|A||Q_{n-1}|^{3}.$$
Again by recursion we show that $$|Q_{n}|=|A|^{\frac{3^{n}-1}{2}}|Q_{0}|^{3^{n}}=|A|^{\frac{3^{n+1}-1}{2}}=|C|^{\frac{s3^{n+1}-s}{2}}.$$
So we obtain for all $n\geq0$ :\begin{eqnarray}\label{a}\bigg|\beta -\displaystyle\frac{P_{n}}{Q_{n}}\bigg|=\displaystyle\frac{1}{|C|^{2s}|Q_{n}|^{\frac{4s-1}{s}}}.\end{eqnarray}
Since $\deg A\geq \deg C$ then $s\geq1$. So $\displaystyle\frac{4s-1}{s}=4-\frac{1}{s}>1+\sqrt{3}$. Hence, if we put $\mu=3, \lambda=|A|, \rho=1/|C|^{2s},\gamma=(4s-1)/s$
then $\gamma>1+ \sqrt{\mu}$ and  following Theorem \ref{t} we conclude that
$\nu(\beta)=4-\displaystyle\frac{1}{s}$.
\end{proofname}

\begin{thm}\label{t3} Let $\beta$ be the irrational solution of the equation $(W_{2})$ such that $|\beta|<1$. Assume that $C$ divides $A$.
Then the continued fraction expansion of $\beta$ is $$[b_{0}, b_{1}, \ldots, b_{n},\ldots]$$
such that $b_{0}=0$,  $b_{1}=A/C$ and for all $n\geq2$:

\begin{eqnarray}\label{Z} b_{n}=\displaystyle(\frac{A}{C})^{3^{n-1}}(C)^{\frac{3^{n-1}+(-1)^{n}}{4}}.\end{eqnarray}
Furthermore, $\nu(\beta)=4$.
\end{thm}
\begin{proofname}. We have $|\beta|<1$ then $b_{0}=0$. Let $\beta_{1}=\beta^{-1}$ then $\beta_{1}$ satisfies the equation
$C\beta_{1}^{4}-A\beta_{1}^{3}-1=0$.
Clearly $[\beta_{1}]=b_{1}=A/C$. So the first partial quotient of $\beta_{1}$ is $b_1=A/C$ and
$\beta_{1}=\displaystyle\frac{A}{C} +\frac{1}{\beta_2}$.
We can easily see that $\beta_{1}$ satisfies  $$\beta_{1}^{3}=\displaystyle\frac{1}{-A+C\beta_{1}}=\frac{\beta_2}{C},$$
then $C\beta_{1}^{3}=\beta_2$. So $Cb_{1}^{3}+\displaystyle\frac{C}{\beta_2^{3}}=\beta_2$. Hence $b_2=Cb_{1}^{3}$ and $\beta_3=\beta_2^{3}/C$.
We apply again the same reasoning and we obtain that $\beta_3=\displaystyle\frac{b_2^{3}}{C} +\displaystyle\frac{1}{C\beta_3^{3}}$, so $b_3=b_2^{3}/C$
and $\beta_{4}=C\beta_3^{3}$. By recurrence on $k$ we prove easily that $\beta_{2k}=C\beta_{2k-1}^{3}$, $\beta_{2k+1}=\beta_{2k}^{3}/C$ and
\begin{eqnarray}\label{ZZ}b_{k}=\left\{%
\begin{array}{ll}
    Cb_{k-1}^{3} & \hbox{if $k$ is even ;} \\
    b_{k-1}^{3}/C  & \hbox{if $k$ is odd.} \\
\end{array}%
\right.    \end{eqnarray}
On the other hand, we have  $b_3=b_2^{3}/C=C^{2}b_{1}^{3^{2}}$ and $b_4=Cb_{3}^{3}=C^{3^{2}-2}b_{1}^{3^{3}}$. We remark that $b_3=C^{\frac{3^{2}-1}{4}}b_{1}^{3^{2}}$
and $b_4=C^{\frac{3^{3}+1}{4}}b_{1}^{3^{3}}$. So by a simple recurrence on $k$ we can prove that $b_k=C^{\frac{3^{k-1}+(-1)^{k}}{4}}b_{1}^{3^{k-1}}$.
Then we deduce that the sequences of partial quotients of $\beta$ is given by: $b_{0}=0$, $b_{1}=A/C$ and for all $n\geq2$:
 $$b_{n}=(A/C)^{3^{n-1}}C^{\frac{3^{n-1}+(-1)^{n}}{4}}.$$
Let $a=\deg A$ and $c=\deg C$. We can compute the approximation exponent of $\beta$ by the formula \eqref{A}:
\begin{eqnarray*}\nu(\beta)&=&2+\displaystyle\limsup \frac{3^{n}(a-c)
-\frac{3^{n}+(-1)^{n+1}}{4}c}{\sum_{k=1}^{n}(3^{k-1}(a-c)-\frac{3^{k-1}+(-1)^{k}}{4}c)}\\&=&2 +2=4.\end{eqnarray*}
\end{proofname}
In the following Theorem, we will give the value of $\nu(\beta)$ for $\beta$ satisfying the equation $(W_{2})$ with the condition on the coefficients of
this equation that is: $C$ does not divides $A$.
\begin{thm}\label{t3} Let $\beta$ be the irrational solution of equation $(W_{2})$ such that $|\beta|<1$. Assume that $C$ does not divide $A$.
Suppose that $|A|=|C|^{s}$ with
$s>\displaystyle\frac{3}{3-\sqrt{3}}$. Then $$\nu(\alpha)=4-\displaystyle\frac{3}{s}$$
\end{thm}
\begin{proofname}. Let $\beta_{1}$ be the first complete quotient of $\beta$. We can easily see that
$\beta_{1}$  satisfies the equation $C\beta_{1}^{4}-A\beta_{1}^{3}-1=0$ and $|\beta_{1}|=|A/C|$. So we have $|\beta_{1}|=|C|^{s-1}$.\\
We consider the following sequence: $P_{0}=A\ , \ Q_{0}=C$ and for $n\geq1$$$P_{n}=AP_{n-1}^{3} +Q_{n-1}^{3}$$$$Q_{n}= CP_{n-1}^{3}.$$
It is easily to see that $\beta_{1}=\displaystyle\frac{1}{C\beta_{1}^{3}} +\frac{A}{C}$ and
$\displaystyle\frac{P_{n}}{Q_{n}}=\frac{Q_{n-1}^{3}}{CP_{n-1}^{3}}+\frac{A}{C}$.
Then for all $n\geq0$:\\
$\bigg|\beta_{1}-\displaystyle\frac{P_{n}}{Q_{n}}\bigg|=\bigg|\frac{1}{C\beta_{1}^{3}}-\frac{Q_{n-1}^{3}}{CP_{n-1}^{3}}\bigg|=
\frac{1}{|C||\beta_{1}|^{3}}\bigg|\frac{\beta_{1}}{|\beta_{1}|}-\frac{P_{n-1}}{|\beta_{1}|Q_{n-1}}\bigg|^{3}
=\frac{1}{|C||\beta_{1}|^{6}}\bigg|\beta_{1}-\frac{P_{n-1}}{Q_{n-1}}\bigg|^{3}.$\\
We show by recursion that for all $n\geq0$:
$$\bigg|\beta_{1}-\displaystyle\frac{P_{n}}{Q_{n}}\bigg|=|C|^{-\frac{(3^{n}-1)}{2}}|\beta_{1}|^{-\frac{6(3^{n}-1)}{2}}\bigg|\beta_{1}-\frac{P_{0}}{Q_{0}}\bigg|^{3^{n}}$$
\noindent since $\bigg|\beta_{1}-\displaystyle\frac{P_{0}}{Q_{0}}\bigg|=\bigg|\beta_{1}-\frac{A}{C}\bigg|=\frac{1}{|C||\beta_{1}|^{3}}$ then
$\bigg|\beta_{1}-\displaystyle\frac{P_{0}}{Q_{0}}\bigg|^{3^{n}}=|C|^{-3^{n}}|\beta_{1}|^{-3^{n+1}}$. So
$$\bigg|\beta_{1}-\displaystyle\frac{P_{n}}{Q_{n}}\bigg|=|C|^{-\frac{3^{n+1}-1}{2}}|\beta_{1}|^{-\frac{3^{n+2}+3^{n+1}-6}{2}}=
|C|^{-\frac{(s-1)3^{n+2}+s3^{n+1}-6(s-1)-1}{2}}.$$
On the other hand, we have for all  $n\geq1$ $Q_{n}= CP_{n-1}^{3}$ and since $|P_{n-1}|=|C|^{s-1}|Q_{n-1}|$ then $$|Q_{n}|=|C|^{3s-2}|Q_{n-1}|^{3}.$$
Again by recursion we show that $$|Q_{n}|=|C|^{\frac{(3s-2)(3^{n}-1)}{2}}|Q_{0}|^{3^{n}}=|C|^{\frac{s3^{n+1}-3s+2}{2}}.$$
So we obtain for all $n\geq0$: $$\bigg|\beta_{1}-\displaystyle\frac{P_{n}}{Q_{n}}\bigg|=|C|^{-\frac{3(s-1)^{2}}{s}}|Q_{n}|^{-\frac{4s-3}{s}}.$$
We can verifies that if $s>\displaystyle\frac{3}{3-\sqrt{3}}$ then $\displaystyle\frac{4s-3}{s}>1+\sqrt{3}$.
Hence by Theorem \eqref{t} we conclude that $\nu(\beta_{1})=4-\displaystyle\frac{3}{s}=4-\displaystyle\frac{3\deg C}{\deg A}$.
\end{proofname}
\section{Diophantine approximation of some not-hyperquadratic power series of degree four over $\mathbb{F}_{3}(T)$}
Now we will give a family of formal power series, defined by their continued fraction expansion, having a minimum value of approximation exponent.
Before this, we recall some usual properties of continued fractions. If $\Omega_{k}=a_{1}, a_{2}, \ldots, a_{k}$ is a sequence of polynomials, we denote $\widetilde{\Omega_{k}}$ the sequence obtained by
reversing the terms of $\Omega_{k}$, i.e, $\widetilde{\Omega_{k}}=a_{k}, a_{k-1}, \ldots, a_{1}$. If $B$ is nonzero element of $\mathbb{F}_{3}[T]$
such that $B$ divides $a_{i}$ for all odd $i$ then  $B^{-1}\Omega_{k}=B^{-1}a_{1}, Ba_{2}, \ldots,$\\
$B^{(-1)^{k}}a_{k}$. Also, if $B$ is nonzero element
of $\mathbb{F}_{3}[T]$ such that $B$ divides $a_{i}$ for all even $i$ then  $B\Omega_{k}=Ba_{1}, B^{-1}a_{2}, \ldots, B^{(-1)^{k-1}}a_{k}$.
In particular, if $\epsilon$ is nonzero element of $\mathbb{F}_{3}$ then we write $\epsilon\Omega_{k}$ for
$\epsilon a_{1}, \epsilon ^{-1}a_{2}, \ldots, \epsilon^{(-1)^{k-1}}a_{k}$. Moreover, in $\mathbb{F}_{3}$ we have $ \epsilon ^{-1}=\epsilon$.
\begin{thm}\label{tf4}Let $A$ and $C$ be two nonzero polynomials in $\mathbb{F}_{3}[T]$ such that
$A$ is not constant, $\deg C\leq \deg A$ and $C$ divides $A$. Let us define the sequence $(\Omega_{n})_{n\geq1}$ of finite sequences of elements of
$\mathbb{F}_{3}[T]$ recursively by
$\Omega_{0}=\emptyset$, $\Omega_{1}=A^{2}$ and for all $n\geq0$
\begin{eqnarray}\label{rrrr}\left\{%
\begin{array}{ll}
    \Omega_{2n +1}= \Omega_{2n}, 2A^{2}, \displaystyle\frac{1}{C^{2}}\Omega_{2n-1}^{(3)}, 2A^{2}, \widetilde{\Omega}_{2n} \\
    \Omega_{2n +2}= \Omega_{2n +1}, A^{2}/C, \displaystyle\frac{1}{2C}\Omega_{2n}^{(3)}, 2A^{2}, \displaystyle\frac{1}{2C}\Omega_{2n +1} \\
\end{array}%
\right.  \end{eqnarray}
Let $\Omega_{\infty}=\displaystyle\lim_{n\longrightarrow \infty}\Omega_{n}$. Let $\theta\in\mathbb{F}_{3}((T^{-1}))$ such that
$\theta=[0, a_{1},\ldots, a_{n},\ldots]=[0,\Omega_{\infty}]$. Then, there
exist explicitly positive numbers $\lambda_{1}$
and $\lambda_{2}$ such that for some rationals $P/Q$ with $|Q|$ arbitrarily large, we have
\begin{eqnarray}\label{zi}|\theta -P/Q|\leq|Q|^{-(2 +\lambda_{1}/\sqrt{\deg Q})}\end{eqnarray}
and, for all rationals $P/Q$ with $|Q|$ sufficiently large, we have
\begin{eqnarray}\label{22}|\theta -P/Q|\geq|Q|^{-(2 +\lambda_{2}/\sqrt{\deg Q})}\end{eqnarray}
where $\lambda_{1}=2/\sqrt{3}$ and $\lambda_{2}>2/\sqrt{3}$.
\end{thm}
\begin{proofname}. We have $\Omega_{2}=A^{2}, A^{2}, 2A^{2}, 2A^{2}/C$. Since $C$ divides $a_{1}=A^{2}$ and $a_{3}=2A^{2}$, then $C$
divides the partial quotient of odd index in $\Omega_{2}$. Suppose
that $C$ divides the partial quotients with odd index in $\Omega_{n}$ for an even $n$. From \eqref{rrrr} we have
$\Omega_{n +1}= \Omega_{n}, 2A^{2}, \displaystyle\frac{1}{C^{2}}\Omega_{n-1}^{(3)}, 2A^{2}, \widetilde{\Omega}_{n}$. As $\Omega_{n}$
has even number of partial quotients then $2A^{2}$ is a  partial quotient with odd index and $C$ divides it. Furthermore,
$\displaystyle\frac{1}{C^{2}}\Omega_{n-1}^{(3)}$
has odd number of partial quotients and begins with a partial quotient with even index, then the  partial quotient $2A^{2}$, coming
after it, has an odd index and $C$ divides it. Finally, as $C$ divides all the partial quotients with odd index in $\Omega_{n}$ then it divides
all partial quotients with even index in $\widetilde{\Omega}_{n}$. So we can compute all the partial quotients of $C^{-1}\Omega_{n+1}$ which is
$$C^{-1}\Omega_{n +1}= C^{-1}\Omega_{n}, 2A^{2}/C, \displaystyle\frac{1}{C}\Omega_{n-1}^{(3)}, 2A^{2}/C, C\widetilde{\Omega}_{n}.$$
By recursion, we prove that we can compute all partial quotients of $C^{-1}\Omega_{n}$ for all $n$.\\
We put $a=\deg A$ and $c=\deg C$. Let us define for each $n\geq 0$, the sequence $\Omega_{n}^{*}$ of the degrees of the elements of $\Omega_{n}$.
The sequence $c^{-1}\Omega_{n}^{*}$ is
the  sequence of degree of $C^{-1}\Omega_{n}$.
We get, from the recursive definition \eqref{rrrr}, $\Omega_{0}^{*}=\emptyset$ and\\
$\Omega_{1}^{*}=2a$\\
$\Omega_{2}^{*}=2a,2a-c,2a,2a-c$\\
$\Omega_{3}^{*}=2a,2a-c,2a,2a-c,2a,(6a-2c),2a,2a-c,2a,2a-c,2a$\\
$\Omega_{4}^{*}=\Omega_{3}^{*},2a-c,(6a-c,6a-2c,6a-c,6a-2c),2a,c^{-1}\Omega_{3}^{*}$\\
$\Omega_{5}^{*}=\Omega_{4}^{*},2a,6a-2c,6a-c,6a-2c,6a-c,6a-2c,(18a-4c),6a-2c,6a-c,6a-2c,$\\
$6a-c,6a-2c,2a,\widetilde{\Omega}_{4}^{*}$\\
From the definition of the approximation exponent, we see that we shall use, for all $k\geq1$, $\Omega^{*}_{2k +1}$
to compute the value of the approximation exponent. Again, from \eqref{rrrr} and by induction on $k$ we see that $\Omega^{*}_{2k +1}$ has an odd number of terms, has
$2(3^{k}a-\displaystyle\frac{3^{k} +(-1)^{k+1}}{4}c)$ as the central term, and is reversible.\\
For $k\geq1$ we put $d_{k}=\deg a_{k}$ and $P/Q=[a_{1},\ldots, a_{k}]$. We define $k_{i}=\inf\{k\geq 1; d_{k}=2(3^{i}a-\displaystyle\frac{3^{i} +(-1)^{i+1}}{4}c)\}$. So we have
\begin{eqnarray}\label{30}\displaystyle\sum_{a_{k}\in \Omega_{2i +1} }d_{k}=\displaystyle2(3^{i}a-\displaystyle\frac{3^{i} +(-1)^{i+1}}{4}c) +2\sum_{k<k_{i}}d_{k}.\end{eqnarray}
Now we put $D_{n}=\displaystyle\sum_{a_{k}\in \Omega_{n}}d_{k}$. Furthermore, we have $D_{n}=\deg\Omega_{n}=2\deg Q_{n}$.
$$D_{2i +1}=2\displaystyle\sum_{k=1}^{2i +1}(3^{k-1}a-\frac{3^{k-1} +(-1)^{k}}{4}c)=2.(\frac{3^{2i +1}-1}{2}a-\frac{3^{2i +1}-3}{8}c).$$
Hence, if $(U_{k}/V_{k})_{k\geq0}$ is the sequence of convergents of $\theta$, the relation \eqref{30} implies, for $i\geq 1$,
\begin{eqnarray*}
\deg V_{k_{i}-1}&=&\displaystyle\sum_{k< k_{i} }d_{k}\\&=
&(D_{2i +1}-2(3^{i}a-\displaystyle\frac{3^{i} +(-1)^{i+1}}{4}c))/2\\&=&(\frac{3^{2i +1}-23^{i}-1}{2})(a-\frac{c}{4})+(1+(-1)^{i})\frac{c}{4}.
\end{eqnarray*}
We can easily verify that $2(3^{i}a-\displaystyle\frac{3^{i} +(-1)^{i+1}}{4}c)\geq 2/\sqrt{3}\sqrt{\deg V_{k_{i}-1}}$, which gives that
$$|T|^{-2(3^{i}a-\displaystyle\frac{3^{i} +(-1)^{i+1}}{4}c)}\leq |V_{k_{i}-1}|^{2/\sqrt{3\deg V_{k_{i}-1}}}.$$
On the other hand, for $i\geq 1$, we have
$$|\theta -U_{k_{i}-1}/V_{k_{i}-1}|=|T|^{-2(3^{i}a-\displaystyle\frac{3^{i} +(-1)^{i+1}}{4}c)}|V_{k_{i}-1}|^{-2}.$$
So, we obtain the desired inequality for $P/Q=U_{k_{i}-1}/V_{k_{i}-1}$ and for $i\geq 1$, with $\lambda_{1}=2/\sqrt{3}$.\\
Furthermore if $U_{k}/V_{k}$ is a convergent to $\theta$, then
$$\deg V_{k_{i}-1}\leq \deg V_{k}< \deg V_{k_{i +1}-1} \ \ implies\ \ |\theta -U_{k}/V_{k}|=|T|^{d_{k +1}}|U_{k}|^{-2}$$
As $\limsup\displaystyle\frac{2(3^{i}a-\displaystyle\frac{3^{i} +(-1)^{i+1}}{4}c)}{\sqrt{\deg V_{k_{i}-1}}}=2/\sqrt{3}$, then, if $\lambda_{2}>2/\sqrt{3}$, we can write
$$2(3^{i}a-\displaystyle\frac{3^{i} +(-1)^{i+1}}{4}c)<\lambda_{2}\sqrt{\deg V_{k_{i}-1}}\leq \lambda_{2}\sqrt{\deg V_{k}}$$ for $i$ large enough.
It follows that \eqref{22} holds for $U_{k}/V_{k}$ with $k$ large enough. Since the convergents are the best rational approximation, this is also true for
all $P/Q$ with $|Q|$ large enough.
\end{proofname}

Let $\beta\in \mathbb{F}_{3}((T^{-1}))$ be the solution of the equation $(W_{1})$ such that $C$ divides $A$. We know
from the Theorem \ref{t3} that the continued fraction expansion of $\beta$ is: $$[b_{0}, b_{1}, \ldots, b_{n},\ldots]$$
such that $b_{0}=0$, $b_{1}=A$ and for all $n\geq2$:
$$b_{n}=(-1)^{n-1}A^{3^{n-1}}C^{-\frac{3^{n-1}+(-1)^{n}}{4}}.$$
In the next part, we will compute the continued fraction expansion and the approximation exponent of
$\alpha=\beta^{2}=[0, A, -A^{3}/C,\ldots,(-1)^{n-1}A^{3^{n-1}}C^{-\frac{3^{n-1}+(-1)^{n}}{4}}, \ldots]^{2}$. Note that, from the equation $(W_{1})$,
$\beta$ satisfies $\beta=(C\beta^{4}+1)/A$. So $\beta^{2}=(C\beta^{4}+1)^{2}/A^{2}$, which gives that $C^{2}\beta^{8}+2C\beta^{4}+1=A^{2}\beta^{2}$.
Then we deduce that $\alpha$ satisfies the equation $C^{2}\alpha^{4}+2C\alpha^{2}-A^{2}\alpha^{2}+1=0$ which is the equation $(E_{1})$.\\
We set $\alpha=[a_{0}, a_{1}, \ldots, a_{n},\ldots].$ Observe that $a_{0}=0$ from the definition of $\alpha$ since $|\beta|<1$. Then we introduce
the usual two sequences of polynomials of $\mathbb{F}_{3}[T]$, defined inductively by
$$U_{0}=0, \ \ \ U_{1}=1,\ \ \ V_{0}=1,\ \ \ V_{1}=a_{1},$$
$$U_{n}=a_{n}U_{n-1} +U_{n-2},\ \ \ V_{n}=a_{n}V_{n-1} +V_{n-2}$$
for $n\geq2$. So $(U_{n}/V_{n})_{n\geq0}$ is the sequence of the convergents to $\alpha$.\\
Now, in order to compute all the partial quotients of $\alpha$, we need to introduce a series of Lemmas.
\begin{lem}
Let $(P_{n}/Q_{n})_{n\geq0}$ be the sequence of convergents
of $\beta$. Then $P_{0}=0$, $Q_{0}=1$, $P_{1}=1$, $Q_{1}=A$ and for all $n\geq1$:
\begin{eqnarray}\label{ref}\left\{%
\begin{array}{ll}
    P_{2n +1}=Q_{2n}^{3} \\
    Q_{2n +1}=AQ_{2n}^{3}-CP_{2n}^{3}\\
\end{array}%
\right.
\ \ and\ \ \left\{%
\begin{array}{ll}
    P_{2n}=-Q_{2n-1}^{3}/C \\
    Q_{2n}=-(A/C)Q_{2n-1}^{3}+P_{2n-1}^{3}\\
\end{array}%
\right.    \end{eqnarray}
\end{lem}
\begin{proofname}. From the equality \eqref{Z} defining the sequence of partial quotients of $\beta$ we can easily check that
$\displaystyle\frac{P_{2}}{Q_{2}}=[0, A, -A^{3}/C]=\frac{-A^{3}/C}{-A^{4}/C+1}$, and
$\displaystyle\frac{P_{3}}{Q_{3}}=[0, A, -A^{3}/C, A^{9}/C^{2}]=\frac{-A^{12}/C^{3}+1}{-A^{13}/C^{3}+A^{9}/C^{2}+1}$. So $P_{2}=-A^{3}/C=-Q_{1}^{3}/C$,
$Q_{2}=-A^{4}/C+1=-(A/C)Q_{1}^{3}+P_{1}^{3}$, $P_{3}=-A^{12}/C^{3}+1=Q_{2}^{3}$ and $Q_{3}=-A^{13}/C^{3}+A^{9}/C^{2}+1=AQ_{2}^{3}-CP_{2}^{3}$. Hence
\eqref{ref} is satisfied for $n=1$. Suppose that \eqref{ref} is satisfied for $n=l>1$. We know that $P_{2l+2}=b_{2l+2}P_{2l+1} +P_{2l}$ and
$Q_{2l+2}=b_{2l+2}Q_{2l+1} +Q_{2l}$. Then
\begin{eqnarray*}
P_{2l+2}&=&(b_{2l+1}^{3}/-C)Q_{2l}^{3} -Q_{2l-1}^{3}/C=(b_{2l+1}^{3}Q_{2l}^{3} +Q_{2l-1}^{3})/-C\\&=&(b_{2l+1}Q_{2l} +Q_{2l-1})^{3}/-C=-Q_{2l+1}^{3}/C,
\end{eqnarray*}
and
\begin{eqnarray*}Q_{2l+2}&=&(b_{2l+1}^{3}/-C)(AQ_{2l}^{3}-CP_{2l}^{3}) +(-(A/C)Q_{2l-1}^{3}+P_{2l-1}^{3})\\&=&-(A/C)(b_{2l+1}^{3}Q_{2l}^{3} +Q_{2l-1}^{3}) +
(b_{2l+1}^{3}P_{2l}^{3} +P_{2l-1}^{3})\\&=&-(A/C)Q_{2l+1}^{3}+P_{2l+1}^{3}.\end{eqnarray*}
So the right part of \eqref{ref} is satisfied for $n=l+1$. Samely, we can obtain the left part.
By induction, we see that \eqref{ref} holds for all $n\ge 1$.
\end{proofname}
We note that the polynomials $P_{n}$ and $Q_{n}$ defined in the previous Lemma will be used throughout the rest
of this section. Also, it is clear that $C$ divides $Q_{n}$ for all $n$ odd integer. Moreover, for the proofs of the following Lemmas, we will
follow \cite{LA2} fairly closely.
\begin{lem}\label{l0} Let $P$ and $Q$ be two polynomials of $\mathbb{F}_{3}[X]$, with $Q\neq0$, and $n$ a positive integer.
Suppose that $PQ_{n}^{2}-QP_{n}^{2}\neq0$. If
\begin{eqnarray}\label{l3}|Q|\leq |Q_{n}|^{2}\ \ \ \ \ and\ \ \ \ |PQ_{n}^{2}-QP_{n}^{2}|<\displaystyle\frac{|Q_{n}|^{2}}{|Q|}\end{eqnarray}
then $P/Q$ is a convergent to $\alpha$. Moreover, if $P$ and $Q$ are coprime and the convergent $P/Q$ is $U_{k}/V_{k}$, then we have
\begin{eqnarray}\label{l4}|a_{k +1}|=|PQ_{n}^{2}-QP_{n}^{2}|^{-1}|Q|^{-1}|Q_{n}|^{2}.\end{eqnarray}
\end{lem}
\begin{proofname}.
We have for $n\geq0$
$$|\beta^{2}-(P_{n}/Q_{n})^{2}|=|\beta-(P_{n}/Q_{n})||\beta+(P_{n}/Q_{n})|.$$
Since $|\beta|=|P_{n}/Q_{n}|=|A|^{-1}$, we have two terms in the sum, each with the absolute value $|A|^{-1}$ and the same dominant coefficient. So this becomes
$$|\beta^{2}-(P_{n}/Q_{n})^{2}|=|\beta-(P_{n}/Q_{n})||A|^{-1}=|Q_{n}Q_{n+1}|^{-1}|A|^{-1}.$$
So \begin{eqnarray}\label{l1}|\beta^{2}-(P_{2n}/Q_{2n})^{2}|=|Q_{2n}Q_{2n+1}|^{-1}|A|^{-1}=|Q_{2n}|^{-4}|A|^{-2},\end{eqnarray}
\begin{eqnarray}\label{l2}|\beta^{2}-(P_{2n+1}/Q_{2n+1})^{2}|=|Q_{2n+1}Q_{2n+2}|^{-1}|A|^{-1}=|Q_{2n+1}|^{-4}|A^{2}/C|^{-1}.\end{eqnarray}
From the equalities \eqref{l1} and \eqref{l2} we have:
$$\displaystyle|\alpha -(P_{n}/Q_{n})^{2}|\leq\frac{1}{|Q_{n}|^{4}|A|}
<\frac{1}{|Q_{n}|^{4}}\leq\frac{1}{|Q_{n}|^{2}|Q|}\leq \frac{|PQ_{n}^{2}-QP_{n}^{2}|}{|Q_{n}|^{2}|Q|}.$$ Hence
$$|\alpha -(P_{n}/Q_{n})^{2}|<|P/Q -(P_{n}/Q_{n})^{2}|.$$
Therefore,
$$|\alpha- P/Q|=|\alpha -(P_{n}/Q_{n})^{2} +(P_{n}/Q_{n})^{2}- P/Q|=|P/Q- (P_{n}/Q_{n})^{2}|$$
and by \eqref{l3}
$$|\alpha- P/Q|<|Q|^{-2}.$$
This shows that $P/Q$ is a convergent to $\alpha$. Now if $P$ and $Q$ are coprime and $P/Q=U_{k}/V_{k}$, we have $|Q|=|V_{k}|$. Besides, we know that
$$|\alpha-U_{k}/V_{k}|=|V_{k}|^{-2}|a_{k +1}|^{-1}.$$
Since $$|\alpha-U_{k}/V_{k}|=|P/Q -(P_{n}/Q_{n})^{2}|,$$
then \eqref{l4}holds and so we obtain the desired result.
\end{proofname}
\noindent We denote by $a=\deg A$ and $c=\deg C$.
\begin{lem}\label{l1} We consider the following sequences of rational functions:

\begin{eqnarray}\label{-21}\displaystyle\frac{R_{1,n}}{S_{1,n}}=P_{n}^{2}/Q_{n}^{2},\ \ for\ \ all \ \ n\geq1.\end{eqnarray}

$$\displaystyle\frac{R_{2,n}}{S_{2,n}}=\left\{%
\begin{array}{ll}
    C^{-1}P_{n}^{2}Q_{n}^{2}/(C^{-1}Q_{n}^{4} +1) & \hbox{if $n$ is odd ;}\\
    P_{n}^{2}Q_{n}^{2}/(Q_{n}^{4} +1) & \hbox{if $n$ is even.}\\
\end{array}%
\right.    $$

$$\displaystyle\frac{R_{3,n}}{S_{3,n}}=\left\{%
\begin{array}{ll}
    P_{n}^{2}(-Q_{n}^{4}C^{-1}+1)/-C^{-1}Q_{n}^{6} & \hbox{if $n$ is odd ;} \\
    P_{n}^{2}(Q_{n}^{4} +2)/Q_{n}^{6} & \hbox{if $n$ is even.} \\
\end{array}%
\right.    $$
Then for all $ 1\leq i\leq3$, $R_{i,n}/S_{i,n}$ is a convergent to $\alpha$. Further $R_{i,n}$ and $S_{i,n}$ are coprime, and if we put $m(i,n)$
the integer
such that $U_{m(i,n)}/V_{m(i,n)}=R_{i,n}/S_{i,n}$ then:
$$\left\{%
\begin{array}{lll}
    \deg a_{m(1,n)+1}=2a\ \  if\ \  n\ \  is\ \  even\ \  and\ \  \deg a_{m(1,n)+1}=2a -c\ \  if\ \  n\ \  is\ \  odd, \\
    \deg a_{m(2,n)+1}=2a\ \  for\ \  all\ \  n, \\
    \deg a_{m(3,n)+1}=2a\ \  if\ \  n\ \  is\ \  even\ \  and\ \  \deg a_{m(3,n)+1}=2a-c\ \  if\ \  n\ \  is\ \  odd.\\
\end{array}%
\right.    $$
Moreover, we have
$\displaystyle\frac{R_{3,n}}{S_{3,n}}$ is the convergent which comes before $\displaystyle\frac{R_{1,n+1}}{S_{1,n+1}}$, i.e
\begin{eqnarray}\label{-1}U_{m(1,n+1)-1}/V_{m(1,n+1)-1}=R_{3,n}/S_{3,n}\end{eqnarray}
\end{lem}
\begin{proofname}.
The equalities \eqref{l1} and \eqref{l2} gives that for all $n\geq1$, $\displaystyle\frac{R_{1,n}}{S_{1,n}}=P_{n}^{2}/Q_{n}^{2}$ is a sequence
convergent of $\alpha$ such that $\deg a_{m(1,n)+1}=2a$ if $n$ is even and $\deg a_{m(1,n)+1}=2a -c$ if $n$ is odd.\\
On the other hand we have $|S_{2,n}|=|C|^{-1}|Q_{n}|^{4}\leq |Q_{n+1}|^{2}$ if $n$ is odd and $|S_{2,n}|=|Q_{n}|^{4}\leq |Q_{n+1}|^{2}$ if $n$ is
even. Moreover, for odd  $n$ we have $|S_{3,n}|=|Q_{n}|^{6}|C|^{-1}\leq |Q_{n+1}|^{2}$, and for
even $n$ we have $|S_{3,n}|=|Q_{n}|^{6}\leq |Q_{n+1}|^{2}$ then we obtain the first part of the condition
\eqref{l3}.\\
\emph{*)For odd $n$:}\\
\begin{eqnarray*}
R_{2,n}Q_{n +1}^{2} -S_{2,n}P_{n +1}^{2}&=&\displaystyle P_{n}^{2}\frac{Q_{n}^{2}}{C}Q_{n +1}^{2}-(\frac{Q_{n}^{4}}{C} +1)P_{n +1}^{2}\\&=
&(1 -P_{n+1}Q_{n})^{2}\frac{Q_{n}^{2}}{C}-(\frac{Q_{n}^{4}}{C} +1)\frac{Q_{n}^{6}}{C^{2}}\\&=&
(1 -\frac{Q_{n}^{4}}{C})^{2}\frac{Q_{n}^{2}}{C}-(\frac{Q_{n}^{4}}{C} +1)\frac{Q_{n}^{6}}{C^{2}}\\&=&\frac{Q_{n}^{2}}{C}.
\end{eqnarray*}
Let $H$ be a common divisor to $R_{2,n}$ and $S_{2,n}$ then $H$ divides $\displaystyle \frac{Q_{n}^{2}}{C}$ and so $H$ divides $\displaystyle \frac{Q_{n}^{4}}{C}$.
Since $H$ divides $S_{2,n}=\displaystyle \frac{Q_{n}^{4}}{C} +1$ then $H$ divides 1. Thus,  $R_{2,n}$ and $S_{2,n}$ are coprime. On the other hand,
$$|R_{2,n}Q_{n +1}^{2} -S_{2,n}P_{n +1}^{2}|=\frac{|Q_{n}|^{2}}{|C|}< \displaystyle\frac{|Q_{n+1}|^{2}}{|S_{2,n}|}=|A/C|^{2}|Q_{n}|^{2}$$
So $\displaystyle\frac{R_{2,n}}{S_{2,n}}$ is a convergent to $\alpha$ and \\
$$|a_{m(2,n)+1}|=|Q_{n}|^{-2}|C|^{2}|Q_{n}|^{-4}|Q_{n+1}|^{2}=|Q_{n}|^{-2}|C|^{2}|Q_{n}|^{-4}|A/C|^{2}|Q_{n}|^{6}=|A|^{2}.$$
\emph{*)For even $n$:}
\begin{eqnarray*}
|R_{2,n}Q_{n +1}^{2} -S_{2,n}P_{n +1}^{2}|&=&|P_{n}^{2}Q_{n}^{2}Q_{n +1}^{2}-(Q_{n}^{4} +1)P_{n +1}^{2}|\\&=
&|(1 +P_{n+1}Q_{n})^{2}Q_{n}^{2}-(Q_{n}^{4} +1)Q_{n}^{6}|\\&=&|(1 +Q_{n}^{4})^{2}Q_{n}^{2}-(Q_{n}^{4} +1)Q_{n}^{6}|\\&=&
|Q_{n}|^{2}< \displaystyle\frac{|Q_{n+1}|^{2}}{|S_{2,n}|}=|A||Q_{n}|^{2}.
\end{eqnarray*}
Moreover, it is clear that $R_{2,n}$ and $S_{2,n}$ are coprime.
So $\displaystyle\frac{R_{2,n}}{S_{2,n}}$ is a convergent to $\alpha$ and \\
$$|a_{m(2,n)+1}|=|Q_{n}|^{-2}|Q_{n}|^{-4}|Q_{n+1}|^{2}=|Q_{n}|^{-2}|Q_{n}|^{-4}|A|^{2}|Q_{n}|^{6}=|A|^{2}.$$
\emph{*)For odd $n$:}
\begin{eqnarray*}
|R_{3,n}Q_{n +1}^{2} -S_{3,n}P_{n +1}^{2}|&=&\bigg|P_{n}^{2}(-\frac{Q_{n}^{4}}{C}+1)Q_{n +1}^{2}+\frac{Q_{n}^{6}}{C}P_{n +1}^{2}\bigg|\\&=
&\bigg|-(1 -P_{n+1}Q_{n})^{2}(\frac{Q_{n}^{4}}{C}-1)+\frac{Q_{n}^{6}}{C}\frac{Q_{n}^{6}}{C^{2}}\bigg|\\&=&
\bigg|-(1 -\frac{Q_{n}^{4}}{C})^{2}(\frac{Q_{n}^{4}}{C}-1)+\frac{Q_{n}^{6}}{C}\frac{Q_{n}^{6}}{C^{2}}\bigg|\\&=&
1\leq \displaystyle\frac{|Q_{n+1}|^{2}}{|S_{3,n}|}=|A^{2}C|.
\end{eqnarray*}
This gives that $R_{3,n}$ and $S_{3,n}$ are coprime. So $\displaystyle\frac{R_{3,n}}{S_{3,n}}$ is a convergent to $\alpha$ and \\
$$|a_{m(3,n)+1}|=|Q_{n}|^{-6}|C||Q_{n+1}|^{2}=|Q_{n}|^{-6}|C||A/C|^{2}|Q_{n}|^{6}=|A^{2}/C|.$$
\emph{*)For even $n$:}
\begin{eqnarray*}
|R_{3,n}Q_{n +1}^{2} -S_{3,n}P_{n +1}^{2}|&=&|P_{n}^{2}(Q_{n}^{4} +2)Q_{n +1}^{2}-Q_{n}^{6}P_{n +1}^{2}|\\&=
&|(1 -P_{n+1}Q_{n})^{2}(Q_{n}^{4} +2)-Q_{n}^{6}Q_{n}^{6}|\\&=&|(1 -Q_{n}^{4})^{2}(Q_{n}^{4} +2)-Q_{n}^{12}|\\&=&
1\leq \displaystyle\frac{|Q_{n+1}|^{2}}{|S_{3,n}|}=|A|^{2}.
\end{eqnarray*}
This gives that $R_{3,n}$ and $S_{3,n}$ are coprime. So $\displaystyle\frac{R_{3,n}}{S_{3,n}}$ is a convergent to $\alpha$ and \\
$$|a_{m(3,n)+1}|=|Q_{n}|^{-6}|Q_{n+1}|^{2}=|Q_{n}|^{-6}|A|^{2}|Q_{n}|^{6}=|A|^{2}.$$\\
Furthermore, we note that we have $|S_{1, n+1}|=|Q_{n+1}|^{2}=|A|^{2}|Q_{n}^{6}|=|S_{3,n}||A|^{2}$ for even $n$ and
$|S_{1, n+1}|=|Q_{n+1}|^{2}=|A/C|^{2}|Q_{n}|^{6}=|S_{3,n}||A^{2}/C|$ for odd $n$, this leads to deduce that
$R_{3,n}/S_{3,n}$ is the convergent  coming before $R_{1,n+1}/S_{1,n +1}$.
\end{proofname}

\noindent We introduce $\Omega_{1, n}$, $\Omega_{2, n}$ and $\Omega_{3, n}$ the sequences of partial quotients which represent
respectively the convergents $\displaystyle\frac{R_{1,n}}{S_{1,n}}$, $\displaystyle\frac{R_{2,n}}{S_{2,n}}$ and $\displaystyle\frac{R_{3,n}}{S_{3,n}}$.
Then:
$$R_{i,n}/S_{i,n}=[0, \Omega_{i, n}]\ \ and \ \ \Omega_{i, n}=a_{1},\ldots, a_{m(i, n)}\ \ for\ \ n\geq 0\ \ and\ \ 1\leq i\leq3.$$
We have $R_{1,1}/S_{1,1}=1/A^{2}$ so $\Omega_{1, 1}=a_{1}=A^{2}$. Further, $$R_{1,2}/S_{1,2}=P_{2}^{2}/Q_{2}^{2}=[0, A^{2},A^{2}/C,2A^{2},2A^{2}/C]$$ so
$\Omega_{1, 2}=A^{2},A^{2}/C,2A^{2},2A^{2}/C$.\\
Note that for $n\geq1$, we have $1<|S_{1,n}|<|S_{2,n}|<|S_{3,n}|$ then $m(1, n)<m(2,n)<m(3,n)$. We put
$\Lambda_{2, n}=a_{m(1,n) +2},\ldots, a_{m(2,n)}$ and $\Lambda_{3, n}=a_{m(2,n) +2},\ldots, a_{m(3,n)}$. Then, from the previous Lemma we can write for $n\geq1$:
$$\Omega_{2, n}=\Omega_{1, n}, a_{m(1,n) +1}, \Lambda_{2, n},$$
$$\Omega_{3, n}=\Omega_{1, n}, a_{m(1,n) +1}, \Lambda_{2, n}, a_{m(2,n) +1}, \Lambda_{3, n}.$$
So, from \eqref{-1} we can write for $n\geq1$:
\begin{eqnarray}\label{-2}
\Omega_{1, n+1}=\Omega_{1, n}, a_{m(1,n) +1}, \Lambda_{2, n}, a_{m(2,n) +1}, \Lambda_{3, n}, a_{m(3,n) +1}.
\end{eqnarray}
On the other hand, observations by computers of the first few hundred of partial quotients of the solution $\alpha$ of $(E_{1})$ show that
 $C$ divides all partial quotients with odd index of any sequence
$\Omega_{k}=a_{1}, a_{2}, \ldots, a_{k}$ and we can compute the sequence of partial quotients of $C^{-1}\Omega_{k}$, as we have describe above.
So, we admit this in the following Lemma, more precisely, equality \eqref{nik} below.
However, we are not able to provide a proof. For this reason, we will state our last result as a conjecture and we will expose this problem as an open question
at the end of this section.
\begin{lem}\label{l2}
There exists, a sequences $(\epsilon_{n})_{n\geq1}$  of nonzero element of $\mathbb{F}_{3}$, such that\\
\emph{-)For even $n$:}
\begin{eqnarray}\label{100}a_{m(1, n)-k}=\epsilon_{n}\displaystyle\frac{1}{(2C)^{(-1)^{k}}}a_{k +1}\end{eqnarray}
for each $(k, n)$ with $0\leq k\leq m(1,n)-1; n\geq1$. Further, we have for $n\geq2$,
\begin{eqnarray*}
\left\{%
\begin{array}{ll}
    \Lambda_{3, n}, a_{m(3,n) +1}=\epsilon_{n+1}\widetilde{\Omega}_{1, n}=\displaystyle\frac{1}{2C}\Omega_{1, n}; &
    \Lambda_{2,n}=\epsilon_{n+1}\widetilde{\Lambda}_{2,n} \\
    a_{m(3,n) +1}=\epsilon_{n+1}A^{2}; & a_{m(1,n) +1}=\epsilon_{n+1}a_{m(2,n) +1}\\
\end{array}%
\right.
\end{eqnarray*}
\emph{-)For odd $n$:}
\begin{eqnarray}\label{102}a_{m(1, n)-k}=\epsilon_{n}a_{k +1}\end{eqnarray}
for each $(k, n)$ with $0\leq k\leq m(1,n)-1; n\geq1$. Further we have for $n\geq2$,
\begin{eqnarray*}
\left\{%
\begin{array}{ll}
    \Lambda_{3, n}, a_{m(3,n) +1}=\epsilon_{n+1}\widetilde{\Omega}_{1, n}=\Omega_{1, n}; &
     \Lambda_{2,n}=\epsilon_{n+1}\displaystyle\frac{1}{2C}\widetilde{\Lambda}_{2,n}\\
    a_{m(3,n) +1}=\epsilon_{n+1}2A^{2}/C; & a_{m(1,n) +1}=\epsilon_{n+1}\displaystyle\frac{1}{2C}a_{m(2,n) +1} \\
\end{array}%
\right.
\end{eqnarray*}
\end{lem}
\begin{proofname}. \emph{If $n$ is even:}
By \eqref{-21} and \eqref{-1}, we can write
\begin{eqnarray}\label{l5}U_{m(1,n)}=\epsilon'_{n}P_{n}^{2},\ \ \ \ \ \ V_{m(1,n)}=\epsilon'_{n}Q_{n}^{2}\end{eqnarray}
and
\begin{eqnarray}\label{l6}U_{m(1,n)-1}=\epsilon''_{n}P_{n-1}^{2}(-Q_{n-1}^{4}C^{-1} +1),\ \ \ \ \ \ V_{m(1,n)-1}=-\epsilon''_{n}Q_{n-1}^{6}/C=\epsilon''_{n}2CP_{n}^{2}\end{eqnarray}
where $\epsilon'_{n}$ and $\epsilon''_{n}$ are nonzero elements of $\mathbb{F}_{3}$. We write $\epsilon_{n}=\epsilon'_{n}/\epsilon''_{n}$.\\
We can write $V_{m(1,n)}/V_{m(1,n)-1}=[a_{m(1,n)},a_{m(1,n)-1},\ldots, a_{1}]$. On the other hand, by \eqref{l5} and \eqref{l6}, we have
$$\displaystyle\frac{V_{m(1,n)}}{V_{m(1,n)-1}}=\epsilon_{n}\frac{1}{2C}\frac{V_{m(1,n)}}{U_{m(1,n)}}=\epsilon_{n}\frac{1}{2C}\frac{1}{[0, a_{1},\ldots,a_{m(1,n)}]};$$
therefore:
\begin{eqnarray*}[a_{m(1,n)},a_{m(1,n)-1},\ldots, a_{1}]=\epsilon_{n}\displaystyle\frac{1}{2C}[a_{1},\ldots,a_{m(1,n)}].\end{eqnarray*}
Admit that
\begin{eqnarray}\label{nik}
\epsilon_{n}\displaystyle\frac{1}{2C}[a_{1},\ldots,a_{m(1,n)}]=\epsilon_{n}[(2C)^{-1}a_{1},..,(2C)^{(-1)^{i}}a_{i},..,(2C)^{(-1)^{m(1,n)}}a_{m(1,n)}].\end{eqnarray}
Then we can write $\widetilde{\Omega}_{1, n}=\epsilon_{n}\displaystyle\frac{1}{2C}\Omega_{1, n}$ and we get equality \eqref{100}.\\
\emph{If $n$ is odd:}
Again by \eqref{-21} and \eqref{-1}, we can write
\begin{eqnarray}\label{l7}U_{m(1,n)}=\epsilon'_{n}P_{n}^{2},\ \ \ \ \ \ V_{m(1,n)}=\epsilon'_{n}Q_{n}^{2}\end{eqnarray}
and
\begin{eqnarray}\label{l8}U_{m(1,n)-1}=\epsilon''_{n}P_{n-1}^{2}(Q_{n}^{4} -1),\ \ \ \ \ \ V_{m(1,n)-1}=\epsilon''_{n}Q_{n-1}^{6}=P_{n}^{2}\end{eqnarray}
Then we obtain
$$\displaystyle\frac{V_{m(1,n)}}{V_{m(1,n)-1}}=\epsilon_{n}\frac{V_{m(1,n)}}{U_{m(1,n)}}=\epsilon_{n}\frac{1}{[0, a_{1},\ldots,a_{m(1,n)}]};$$
therefore:
$$[a_{m(1,n)},a_{m(1,n)-1},\ldots, a_{1}]=\epsilon_{n}[ a_{1},\ldots,a_{m(1,n)}]$$
Then we can write $\widetilde{\Omega}_{1, n}=\epsilon_{n}\Omega_{1, n}$ and we get equality \eqref{102}.\end{proofname}
For $n\geq 1$, we put $\Omega_{1, n}=a_{1}, \Lambda_{1, n}$. Hence, for $n\geq 1$, \eqref{-2} becomes
\begin{eqnarray}\label{-3}
\Omega_{1, n+1}=A^{2},\Lambda_{1, n}, a_{m(1,n) +1}, \Lambda_{2, n}, a_{m(2,n) +1}, \Lambda_{3, n}, a_{m(3, n)+1}.
\end{eqnarray}
For each finite sequence of nonzero polynomials, we define its degree as being the sum of the degrees of its terms. We have
$\deg \Omega_{1, n}=\deg S_{1,n}=2\deg Q_{n}$.\\
\emph{*)}If $n$ is even then $\widetilde{\Omega}_{1, n+1}=\epsilon_{n+1}\Omega_{1, n+1}$. Further we have $\deg S_{2,n}=4\deg Q_{n}$ and
$\deg S_{3,n}=6\deg Q_{n}$. As $\deg Q_{n}=3\deg Q_{n-1} +a-c$ then if we put $w_{n}=6\deg Q_{n-1}-2c$ then $\deg \Omega_{1, n}=w_{n}+2a$,
$\deg S_{2,n}=12\deg Q_{n-1}+4a-4c=2a+w_{n}+2a+w_{n}$ and $\deg S_{3,n}=18\deg Q_{n-1}+6a-6c=2a+w_{n}+2a+w_{n}+2a+w_{n}$. As
$\deg a_{m(1,n) +1}=\deg a_{m(2,n) +1}=\deg a_{m(3, n)+1}=2a$ then if we write the sequence of the degrees of the
components in the right side of \eqref{-3}, we obtain the sequence, of 7 terms:$2a,w_{n},2a,w_{n},2a,w_{n},2a$. As this sequence is reversible and
$\widetilde{\Omega}_{1, n+1}=\epsilon_{n+1}\Omega_{1, n+1}$, it is clear that $\widetilde{\Lambda}_{3,n}=\epsilon_{n+1}\Lambda_{1, n}$,
$\Lambda_{2,n}=\epsilon_{n+1}\widetilde{\Lambda}_{2,n}$, $a_{m(3, n)+1}=\epsilon_{n+1}A^{2}$, $a_{m(1,n) +1}=\epsilon_{n+1}a_{m(2,n) +1}$.\\
\emph{*)}If $n$ is odd then $\widetilde{\Omega}_{1, n+1}=\epsilon_{n+1}\displaystyle\frac{1}{2C}\Omega_{1, n+1}$. Further we have
$\deg S_{2,n}=4\deg Q_{n}-c$ and
$\deg S_{3,n}=6\deg Q_{n}-c$. As $\deg Q_{n}=3\deg Q_{n-1}+a$ then if we put $w_{n}=6\deg Q_{n-1}$ then $\deg \Omega_{1, n}=w_{n}+2a$,
$\deg S_{2,n}=12\deg Q_{n-1}+4a-c=2a+w_{n}+2a-c+w_{n}$ and $\deg S_{3,n}=18\deg Q_{n-1}+6a-c=2a+w_{n}+2a-c+w_{n}+2a+w_{n}$. As
$\deg a_{m(1,n) +1}=2a-c$, $\deg a_{m(2,n) +1}=2a$ and $\deg a_{m(3, n)+1}=2a-c$, then if we write the sequence of the degrees of the
components in the right side of \eqref{-3}, we obtain the sequence, of 7 terms:$2a,w_{n},2a-c,w_{n},2a,w_{n},2a-c$. As
$\widetilde{\Omega}_{1, n+1}=\epsilon_{n+1}\displaystyle\frac{1}{2C}\Omega_{1, n+1}$, it is clear that
$\widetilde{\Lambda}_{3,n}=\epsilon_{n+1}2C\Lambda_{1, n}$,
$\Lambda_{2,n}=\epsilon_{n+1}2C\widetilde{\Lambda}_{2,n}$, $a_{m(3, n)+1}=\epsilon_{n+1}A^{2}/2C$,
$a_{m(1,n) +1}=\epsilon_{n+1}\displaystyle\frac{1}{2C}a_{m(2,n) +1}$.

\begin{lem}\label{l9}

There exists, a sequences $(\epsilon_{n})_{n\geq1}$  of nonzero element of $\mathbb{F}_{3}$, such that:\\
\emph{-)For even $n$:} we have
$$\Lambda_{2, n}=(\epsilon_{n}/C^{2})\Omega_{n-1}^{(3)}\ \ and \ \ a_{m(1,n) +1}=\epsilon_{n}2A^{2}.$$
\emph{-)For odd $n$:} we have
$$\Lambda_{2, n}=(\epsilon_{n}/2C)\Omega_{n-1}^{(3)}\ \ and \ \ a_{m(1,n) +1}=\epsilon_{n}2A^{2}/C.$$
\end{lem}
\begin{proofname}.

\emph{*)If $n$ is even:}\\
We have $U_{m(1,n)}/V_{m(1,n)}=[0, \Omega_{1, n}]$, $U_{m(1,n) +1}/V_{m(1,n)+1}=[0, \Omega_{1, n}, a_{m(1,n) +1}]$ and\\
$U_{m(2,n)}/V_{m(2,n)}=[0, \Omega_{1, n}, a_{m(1,n) +1}, \Lambda_{2, n} ]$. If we put $x_{2, n}$, the element of $\mathbb{F}_{3}(T)$ defined by
$[\Lambda_{2, n} ]$, then we have
\begin{eqnarray}\label{26}
\displaystyle\frac{U_{m(2,n)}}{V_{m(2,n)}}=\frac{x_{2, n}U_{m(1,n) +1}+ U_{m(1,n)}}{x_{2, n}V_{m(1,n) +1}+ V_{m(1,n)}}.
\end{eqnarray}
We know that $U_{m(2, n)}/ V_{m(2, n)}=R_{2, n}/ S_{2, n}= P_{n}^{2}Q_{n}^{2}/(Q_{n}^{4} +1)$. So if we put
\begin{eqnarray}\label{26''}
P'=P_{n}^{2}Q_{n}^{2}\ \ \ \ \ and\ \ \ \ \ Q'=Q_{n}^{4} +1
\end{eqnarray}
the equality \eqref{26} gives that:
\begin{eqnarray}\label{26'}
x_{2,n}=\epsilon'_{n}\displaystyle\frac{P_{n}^{2}Q'-Q_{n}^{2}P'}{V_{m(1,n) +1}P'- U_{m(1,n)+1}Q'}.
\end{eqnarray}
We should determine $U_{m(1,n)+1}/V_{m(1,n)+1}$. We use the fact that $R_{3, n-1}/S_{3, n-1}$ and\\
$R_{(1, n)}/S_{(1, n)}$ are, from Lemma \ref{l1},
the two reduced precedes it.\\
Hence we consider the polynomials $P$ and $Q$ of $\mathbb{F}_{3}[T]$, defined by:
\begin{eqnarray}\label{27}
P=2A^{2}P_{n}^{2} +P_{n-1}^{3}Q_{n} \ \ \ and\ \ \ \ Q=2A^{2}Q_{n}^{2} -CP_{n}^{2}.
\end{eqnarray}
We will apply Lemma \ref{l0}, to prove that $P/Q$ is a convergent to $\alpha$. First we have $\deg Q=2\deg Q_{n} +2a$ and then $Q\neq 0$.
From \eqref{27} and \eqref{ref}, we have
$PQ_{n}^{2}-QP_{n}^{2}=P_{n-1}^{3}Q_{n}^{3}+CP_{n}^{4}=P_{n-1}^{3}Q_{n}^{3}-P_{n}^{3}Q_{n-1}^{3}=1$, hence $\gcd(P, Q)=1$.
Since $2\deg Q_{n} +2a\leq 2\deg Q_{n +1}$ for $n\geq 2$, the first part of condition; that is $|Q|< |Q_{n +1}|^{2}$, is satisfied. We should prove that
$|PQ_{n +1}^{2}-QP_{n+1}^{2}|<|Q_{n +1}|^{2}/|Q|$. We put
$$X_{1}=Q_{n +1}^{2}P_{n}^{2} -Q_{n}^{2}P_{n+1}^{2}\ \ \ \ and\ \ \ \ X_{2}=P_{n-1}^{3}Q_{n}Q_{n +1}^{2}+CP_{n}^{2}P_{n+1}^{2}.$$
From \eqref{27}, we have $PQ_{n +1}^{2}- QP_{n +1}^{2}=2A^{2}X_{1} +X_{2}$. Since $P_{n+1}Q_{n}-Q_{n +1}P_{n}=1$, and by \eqref{ref}, we have
$$X_{1}=2Q_{n}P_{n+1} +1=Q_{n}^{4} +1$$
then\\
$X_{2}=Q_{n +1}^{2}P_{n-1}^{3}Q_{n}+CP_{n+1}^{2}P_{n}^{2}=(Q_{n +1}/Q_{n})^{2}(1 -CP_{n}^{4})+CP_{n+1}^{2}P_{n}^{2}$\\
$X_{2}=(Q_{n +1}/Q_{n})^{2} -C(P_{n}/Q_{n})^{2}X_{1}$\\
$X_{2}=(Q_{n +1}/Q_{n})^{2}-C(P_{n}Q_{n})^{2} -C(P_{n}/Q_{n})^{2}.$\\
We put $X=PQ_{n +1}^{2}-QP_{n+1}^{2}$. Since $X=2A^{2}X_{1} +X_{2}$, we have\\
$X=2A^{2} +2A^{2}Q_{n}^{4} +(Q_{n +1}/Q_{n})^{2} -C(P_{n}Q_{n})^{2} -C(P_{n}/Q_{n})^{2}$\\
$X=2A^{2} +2A^{2}Q_{n}^{4} +(AQ_{n}^{2} -CP_{n}^{3}/Q_{n})^{2} -C(P_{n}/Q_{n})^{2}(Q_{n}^{4} +1)$\\
Since $Q_{n}^{4}-AP_{n}Q_{n}^{3}+CP_{n}^{4}=P_{n+1}Q_{n}-Q_{n +1}P_{n}=1$ then\\
$X=2A^{2} +(ACQ_{n}P_{n}^{3} +C^{2}P_{n}^{6}/Q_{n}^{2} -C(P_{n}/Q_{n})^{2}(2Q_{n}^{4}-AP_{n}Q_{n}^{3}+CP_{n}^{4}))$\\
$X-2A^{2}=2ACQ_{n}P_{n}^{3}+CP_{n}^{2}Q_{n}^{2}=CP_{n}^{2}Q_{n}P_{n-1}^{3}.$
Since, for $n\geq2$, $|P_{n-1}^{3}|<|Q_{n}|$ and $|C||P_{n}|^{2}<|Q_{n}|^{2}$, this equality implies:
$$|X|<|Q_{n}^{4}|=\displaystyle\frac{|Q_{n+1}|^{2}}{|Q|}$$
Consequently, $P/Q$ is a convergent to $\alpha$, and since $\deg Q=\deg V_{m(1,n)} +2a$, then it is next $U_{m(1,n)}/V_{m(1,n)}$. We can write
\begin{eqnarray}\label{28}U_{m(1,n)+1}=\eta_{n}P\ \ \ and\ \ \ \ V_{m(1,n)+1}=\eta_{n}Q.\end{eqnarray}
By \eqref{l5}, \eqref{l6} and \eqref{27}, and $\epsilon^{-1}=\epsilon$ for $\epsilon\in \mathbb{F}_{3}$, the first equality of \eqref{28} can be written
$$a_{m(1, n)+1}U_{m(1, n)} +U_{m(1, n)-1}=\eta_{n}\epsilon'_{n}2A^{2}U_{m(1, n)} +\eta_{n}\epsilon''_{n}U_{m(1, n)-1}.$$
Since we have $\deg U_{m(1,n)}> \deg U_{m(1, n)-1}$, it follows that $a_{m(1, n)+1}=\eta_{n}\epsilon'_{n}2A^{2}$ and $\eta_{n}\epsilon''_{n}=1$, i.e
$\eta_{n}=\epsilon''_{n}$. Thus, since $\epsilon_{n}=\epsilon'_{n}\epsilon''_{n}$, we obtain:
$$a_{m(1, n)+1}=\epsilon_{n}2A^{2}.$$

So the equality \eqref{26'} becomes:
\begin{eqnarray}\label{29}
x_{2,n}=\epsilon_{n}\displaystyle\frac{P_{n}^{2}Q'-Q_{n}^{2}P'}{Q P'- P Q'}.
\end{eqnarray}
We are able to compute $x_{2,n}$.\\
$P_{n}^{2}Q'-Q_{n}^{2}P'=P_{n}^{2}(Q_{n}^{4} +1)-Q_{n}^{2}(P_{n}Q_{n})^{2}=P_{n}^{2}=Q_{n-1}^{6}/C^{2}.$\\
From we have\\
$Q P'- P Q'=P_{n}^{2}Q_{n}^{2}(2A^{2}Q_{n}^{2} -CP_{n}^{2})-(Q_{n}^{4} +1)(2A^{2}P_{n}^{2} +P_{n-1}^{3}Q_{n})$\\
$Q P'- P Q'=-CP_{n}^{4}Q_{n}^{2}-Q_{n}^{5}P_{n-1}^{3}-(2A^{2}P_{n}^{2} +P_{n-1}^{3}Q_{n})$\\
$Q P'- P Q'=Q_{n}^{2}(P_{n}Q_{n-1}-Q_{n}P_{n-1})^{3}-(2A^{2}P_{n}^{2}+Q_{n}^{2} -AP_{n}Q_{n})$\\
$Q P'- P Q'=A^{2}P_{n}^{2}+Q_{n}^{2} +AP_{n}Q_{n}$\\
$Q P'- P Q'=(Q_{n}-AP_{n})^{2}=P_{n-1}^{6}$.\\
So \eqref{29} gives that
$$x_{2,n}=\epsilon_{n}/C^{2}(Q_{n -1}/P_{n-1})^{6}.$$
Furthermore
$$[a_{1}, \ldots, a_{m(1,n-1)}]=(Q_{n -1}/P_{n-1})^{2}$$
then
$$(\epsilon_{n}/C^{2})(Q_{n -1}/P_{n-1})^{6}=\epsilon_{n}/C^{2}[a_{1}^{3}, \ldots, a_{m(1,n-1)}^{3}].$$
Thus, we conclude that \eqref{-3}can be written as
\begin{eqnarray}\label{-3'}\Omega_{1, n+1}=\Omega_{1, n}, \epsilon_{n}2A^{2}, (\epsilon_{n}/C^{2})\Omega_{1, n-1}^{(3)},
\epsilon_{n+1}\epsilon_{n}2A^{2}, \epsilon_{n+1}\widetilde{\Omega}_{1, n}.\end{eqnarray}\\
\emph{*)If $n$ is odd:}\\
We have $U_{m(1,n)}/V_{m(1,n)}=[0, \Omega_{1, n}]$, $U_{m(1,n) +1}/V_{m(1,n)+1}=[0, \Omega_{1, n}, a_{m(1,n) +1}]$ and\\
$U_{m(2,n)}/V_{m(2,n)}=[0, \Omega_{1, n}, a_{m(1,n) +1}, \Lambda_{2, n} ]$. If we put $x_{2, n}$, the element of $\mathbb{F}_{3}(T)$ defined by
$[\Lambda_{2, n} ]$, then we have
\begin{eqnarray}\label{31}
\displaystyle\frac{U_{m(2,n)}}{V_{m(2,n)}}=\frac{x_{2, n}U_{m(1,n) +1}+ U_{m(1,n)}}{x_{2, n}V_{m(1,n) +1}+ V_{m(1,n)}}.
\end{eqnarray}
We know that $U_{m(2, n)}/ V_{m(2, n)}=R_{2, n}/ S_{2, n}= C^{-1}P_{n}^{2}Q_{n}^{2}/(C^{-1}Q_{n}^{4} +1)$. So if we put
\begin{eqnarray}\label{31''}
P'=C^{-1}P_{n}^{2}Q_{n}^{2}\ \ \ \ \ and\ \ \ \ \ Q'=C^{-1}Q_{n}^{4} +1
\end{eqnarray}
the equality \eqref{31} gives that:
\begin{eqnarray}\label{31'}
x_{2,n}=\epsilon'_{n}\displaystyle\frac{P_{n}^{2}Q'-Q_{n}^{2}P'}{V_{m(1,n) +1}P'- U_{m(1,n)+1}Q'}.
\end{eqnarray}
We should determine $U_{m(1,n)+1}/V_{m(1,n)+1}$. We use the fact that $R_{3, n-1}/S_{3, n-1}$ and\\ $R_{(1, n)}/S_{(1, n)}$ are, from Lemma \ref{l1} , the two reduced
precedes it.\\
Hence we consider the polynomials $P$ and $Q$ of $\mathbb{F}_{3}[T]$, defined by:
\begin{eqnarray}\label{32}
P=(A^{2}/C)P_{n}^{2} +P_{n-1}^{3}Q_{n} \ \ \ and\ \ \ \ Q=(A^{2}/C)Q_{n}^{2} +P_{n}^{2}.
\end{eqnarray}
We will apply Lemma \ref{l0} to prove that $P/Q$ is a convergent to $\alpha$. First we have $\deg Q=2\deg Q_{n} +2a-c$ and then $Q\neq 0$.
From \eqref{32} and \eqref{ref}, we have
$PQ_{n}^{2}-QP_{n}^{2}=P_{n-1}^{3}Q_{n}^{3}-P_{n}^{4}=P_{n-1}^{3}Q_{n}^{3}-P_{n}^{3}Q_{n-1}^{3}=-1$, hence $\gcd(P, Q)=1$.
Since $2\deg Q_{n} +2a-c\leq 2\deg Q_{n +1}$ for $n\geq 2$, the first part of condition; that is $|Q|< |Q_{n +1}|^{2}$, is satisfied. We should prove that
$|PQ_{n +1}^{2}-QP_{n+1}^{2}|<|Q_{n +1}|^{2}/|Q|$. We put
$$X_{1}=Q_{n +1}^{2}P_{n}^{2} -Q_{n}^{2}P_{n+1}^{2}\ \ \ \ and\ \ \ \ X_{2}=P_{n-1}^{3}Q_{n}Q_{n +1}^{2}-P_{n}^{2}P_{n+1}^{2}.$$
From \eqref{27}, we have $PQ_{n +1}^{2}- QP_{n +1}^{2}=(A^{2}/C)X_{1} +X_{2}$. Since $P_{n+1}Q_{n}-Q_{n +1}P_{n}=-1$, and by \eqref{ref}, we have
$$X_{1}=-Q_{n}P_{n+1} +1=C^{-1}Q_{n}^{4} +1$$
then\\
$X_{2}=Q_{n +1}^{2}P_{n-1}^{3}Q_{n}-P_{n+1}^{2}P_{n}^{2}=(Q_{n +1}/Q_{n})^{2}(-1 +P_{n}^{4})-P_{n+1}^{2}P_{n}^{2}$\\
$X_{2}=-(Q_{n +1}/Q_{n})^{2} +(P_{n}/Q_{n})^{2}X_{1}$\\
$X_{2}=-(Q_{n +1}/Q_{n})^{2}+C^{-1}(P_{n}Q_{n})^{2} +(P_{n}/Q_{n})^{2}.$\\
We put $X=PQ_{n +1}^{2}-QP_{n+1}^{2}$. Since $X=A^{2}/CX_{1} +X_{2}$, we have\\
$X=A^{2}C^{-1} +A^{2}C^{-2}Q_{n}^{4} -(Q_{n +1}/Q_{n})^{2} +C^{-1}(P_{n}Q_{n})^{2} +(P_{n}/Q_{n})^{2}$\\
$X=A^{2}C^{-1} +A^{2}C^{-2}Q_{n}^{4} -(-AC^{-1}Q_{n}^{2} +P_{n}^{3}/Q_{n})^{2} +(P_{n}/Q_{n})^{2}(Q_{n}^{4}/C +1)$\\
Since $2Q_{n}^{4}C^{-1}+AC^{-1}P_{n}Q_{n}^{3}-P_{n}^{4}=P_{n+1}Q_{n}-Q_{n +1}P_{n}=-1$ then\\
$X=A^{2}C^{-1} +2AC^{-1}Q_{n}P_{n}^{3} -P_{n}^{6}/Q_{n}^{2} +(P_{n}/Q_{n})^{2})(-Q_{n}^{4}C^{-1}-AC^{-1}P_{n}Q_{n}^{3}+P_{n}^{4})$\\
$X-A^{2}C^{-1}=AC^{-1}Q_{n}P_{n}^{3}-P_{n}^{2}Q_{n}^{2}C^{-1}=P_{n}^{2}Q_{n}(AC^{-1}P_{n}-C^{-1}Q_{n})=P_{n}^{2}Q_{n}P_{n-1}^{3}.$
Since, for $n\geq2$, $|P_{n-1}|^{3}<|Q_{n}|$ and $|P_{n}|^{2}<|Q_{n}|^{2}/C$, this equality implies:
$$|X|<|Q_{n}^{4}||C|^{-1}=\displaystyle\frac{|Q_{n+1}|^{2}}{|Q|}$$
Consequently, $P/Q$ is a convergent to $\alpha$, and since $\deg Q=\deg V_{m(1,n)} +2a-c$, then it is next $U_{m(1,n)}/V_{m(1,n)}$. We can write
\begin{eqnarray}\label{30}U_{m(1,n)+1}=\eta_{n}P\ \ \ and\ \ \ \ V_{m(1,n)+1}=\eta_{n}Q.\end{eqnarray}

By \eqref{l7}, \eqref{l8} and \eqref{30}, and $\epsilon^{-1}=\epsilon$ for $\epsilon\in \mathbb{F}_{3}$, the first equality of \eqref{30} can be written
$$a_{m(1, n)+1}U_{m(1, n)} +U_{m(1, n)-1}=\eta_{n}\epsilon'_{n}A^{2}C^{-1}U_{m(1, n)} +\eta_{n}\epsilon''_{n}U_{m(1, n)-1}.$$
Since we have $\deg U_{m(1,n)}> \deg U_{m(1, n)-1}$, it follows that $a_{m(1, n)+1}=\eta_{n}\epsilon'_{n}A^{2}C^{-1}$ and $\eta_{n}\epsilon''_{n}=1$, i.e
$\eta_{n}=\epsilon''_{n}$. Thus, since $\epsilon_{n}=\epsilon'_{n}\epsilon''_{n}$, we obtain
$$a_{m(1, n)+1}=\epsilon_{n}A^{2}C^{-1}.$$
So the equality \eqref{31'} becomes:
\begin{eqnarray}\label{32}
x_{2,n}=\epsilon_{n}\displaystyle\frac{P_{n}^{2}Q'-Q_{n}^{2}P'}{Q P'- P Q'}.
\end{eqnarray}
We are able to compute $x_{2,n}$.\\
$P_{n}^{2}Q'-Q_{n}^{2}P'=P_{n}^{2}(C^{-1}Q_{n}^{4} +1)-Q_{n}^{2}C^{-1}(P_{n}Q_{n})^{2}=P_{n}^{2}=Q_{n-1}^{6}.$\\
From we have\\
$Q P'- P Q'=C^{-1}P_{n}^{2}Q_{n}^{2}(A^{2}C^{-1}Q_{n}^{2} +P_{n}^{2})-(C^{-1}Q_{n}^{4} +1)(A^{2}C^{-1}P_{n}^{2} +P_{n-1}^{3}Q_{n})$\\
$Q P'- P Q'=C^{-1}P_{n}^{4}Q_{n}^{2}-C^{-1}Q_{n}^{5}P_{n-1}^{3}-(A^{2}C^{-1}P_{n}^{2} +P_{n-1}^{3}Q_{n})$\\
$Q P'- P Q'=C^{-1}Q_{n}^{2}(P_{n}Q_{n-1}-Q_{n}P_{n-1})^{3}-(A^{2}C^{-1}P_{n}^{2}-C^{-1}Q_{n}^{2} +AP_{n}Q_{n})$\\
$Q P'- P Q'=2A^{2}C^{-1}P_{n}^{2}+2C^{-1}Q_{n}^{2} +2AP_{n}Q_{n}$\\
$Q P'- P Q'=2C(C^{-1}Q_{n}-AC^{-1}P_{n})^{2}=2CP_{n-1}^{6}$.\\
So \eqref{32} gives that
$$x_{2,n}=\epsilon_{n}/2C(Q_{n -1}/P_{n-1})^{6}.$$
Furthermore
$$[a_{1}, \ldots, a_{m(1,n-1)}]=(Q_{n -1}/P_{n-1})^{2}$$
then
$$\epsilon_{n}/2C(Q_{n -1}/P_{n-1})^{6}=(\epsilon_{n}/2C)[a_{1}^{3}, \ldots, a_{m(1,n-1)}^{3}].$$
Thus, we conclude that we can write \eqref{-3} as:
\begin{eqnarray}\label{-3"}\end{eqnarray}$$\Omega_{1, n+1}=\Omega_{1, n}, \epsilon_{n}A^{2}C^{-1}, (\epsilon_{n}/2C)\Omega_{1, n-1}^{(3)},
\epsilon_{n+1}\epsilon_{n}A^{2}, (\epsilon_{n+1}/2C)\widetilde{\Omega}_{1, n}.$$
\end{proofname}
Finally, we have to determine $\epsilon_{n}$ for all $n\geq1$. By Lemmas \ref{l1} and \eqref{102} we have simultaneously $ \epsilon_{n}a^{3}_{m(1, n-1)}=\epsilon_{n+1}\epsilon_{n}a^{3}_{1}$, which implies
$a_{m(1, n-1)}=\epsilon_{n+1}a_{1}$ and $a_{m(1, n-1)}=\epsilon_{n-1}a_{1}$. Therefore, $\epsilon_{n+1}=\epsilon_{n-1}$ for even $n$. We can
verify that we have also $\epsilon_{n+1}=\epsilon_{n-1}$ for odd $n$. Since $\Omega_{1, 2}=1/2C\widetilde{\Omega}_{1, 2}$ and
$\Omega_{1, 3}=\widetilde{\Omega}_{1, 3}$ then $\epsilon_{2}=\epsilon_{3}=1$. So, we obtain $\epsilon_{n}=1$ for all $n\geq1$.
Finally, by Lemma \ref{l2}, the sequence $\Omega_{1, n}$ is
reversible for all $n$ odd, and so $\widetilde{\Omega}_{1,n}=\Omega_{1,n}$. The equality \eqref{-3"} becomes:
$$\Omega_{1, n+1}=\Omega_{1, n}, A^{2}/C, (1/2C)\Omega_{1, n-1}^{(3)},
2A^{2}, (1/2C)\widetilde{\Omega}_{1, n}.$$
The equality \eqref{-3'} becomes
$$\Omega_{1, n+1}=\Omega_{1, n}, 2A^{2}, (1/C^{2})\Omega_{1, n-1}^{(3)},
2A^{2}, \widetilde{\Omega}_{1, n}.$$
So we can deduce the following result.
\begin{conj}\label{tf} Let $\alpha \in \mathbb{F}_{3}((T^{-1}))$ be the formal power series, of strictly negative degree, satisfying $(E_{1})$.
Let $(\Omega_{n})_{n\geq0}$ be a finite sequence of elements of $\mathbb{F}_{3}[T]$, defined by $\Omega_{0}=\emptyset$, $\Omega_{1}=A^{2}$ and for all $n\geq0$:
\begin{eqnarray}\label{rbj}\left\{%
\begin{array}{ll}
    \Omega_{2n +1}= \Omega_{2n}, 2A^{2}, \displaystyle\frac{1}{C^{2}}\Omega_{2n-1}^{(3)}, 2A^{2}, \widetilde{\Omega}_{2n} \\
    \Omega_{2n +2}= \Omega_{2n +1}, A^{2}/C, \displaystyle\frac{1}{2C}\Omega_{2n}^{(3)}, 2A^{2}, \displaystyle\frac{1}{2C}\Omega_{2n +1} \\
\end{array}%
\right.  \end{eqnarray}
Let $\Omega_{\infty}$ be the infinite sequence beginning by $\Omega_{n}$ for all $n\geq1$. Then, the continued fraction expansion of
$\alpha$ is $\alpha=[0, \Omega_{\infty}]$.
\end{conj}
We see that the equality \eqref{rbj} has the same shape as the equality \eqref{rrrr}. So this gives that the formal power series described in Theorem \ref{tf} is
not other than the unique solution of the quartic equation $(E_{1})$.\\
\begin{exm}
Let $(\Omega_{n})_{n\geq1}$ be a finite sequence of elements of $\mathbb{F}_{3}[T]$, defined by $\Omega_{1}=T^{2}$, $\Omega_{2}=T^{2},T,2T^{2},2T$
and for all $n\geq0$
$$\left\{%
\begin{array}{ll}
    \Omega_{2n +1}= \Omega_{2n}, 2T^{2}, \displaystyle\frac{1}{T^{2}}\Omega_{2n-1}^{(3)}, 2T^{2}, \widetilde{\Omega}_{2n} \\
    \Omega_{2n +2}= \Omega_{2n +1}, T, \displaystyle\frac{1}{2T}\Omega_{2n}^{(3)}, 2T^{2}, \displaystyle\frac{1}{2T}\Omega_{2n +1} \\
\end{array}%
\right.    $$
Then, we have from conjecture \eqref{tf}: $[0, \Omega_{\infty}]=[0, T, -T^{2}, T^{7}, \ldots, (-1)^{n-1}T^{\frac{3^{n}+(-1)^{n +1}}{4}},\ldots]^{2}$.

\end{exm}

In fact, the power series $\alpha \in \mathbb{F}_{3}((T^{-1}))$ having $[0, \Omega_{\infty}]$ as continued fraction expansion is the solution of the
equation $(E_{1})$ with $A=C=T$. The partial quotients of $\Omega_{5}$ are:
$$
[T^{2}, T, 2T^{2}, 2T, 2T^{2}, T^{4}, 2T^{2}, 2T, 2T^{2}, T, T^{2}, T,
2T^{5}, 2T^{4}, T^{5}, T^{4}, 2T^{2}, 2T, 2T^{2}, T, T^{2}, T, 2T^{5},
T, $$$$T^{2},T, 2T^{2}, 2T, 2T^{2}, T^{4}, T^{5}, 2T^{4}, 2T^{5},
2T^{4}, T^{14}, 2T^{4}, 2T^{5}, 2T^{4}, T^{5}, T^{4}, 2T^{2}, 2T, 2T^{2},
T, T^{2}, T,$$$$ 2T^{5}, T, T^{2}, T, 2T^{2}, 2T, 2T^{2}, T^{4}, T^{5},
2T^{4}, 2T^{5}, T, T^{2}, T, 2T^{2}, 2T, 2T^{2}, T^{4}, 2T^{2}, 2T,
2T^{2}, T, T^{2}].
$$\\
\noindent Note that in this case $A^{2}=T^{2}, A^{2}/C=T$ and we have:
$$R_{1,1}/S_{1,1}=[0, T^{2}]=[0, \Omega_{1}],\ \ \ \ \ \Omega_{1}=T^{2}.$$
$$\ R_{2,1}/S_{2,1}=[0, T^{2}, T]; \ R_{3,1}/S_{3,1}=[0, T^{2}, T, 2T^{2}]$$
$$R_{1,2}/S_{1,2}=[0, T^{2}, T, 2T^{2}, 2T]=[0, \Omega_{2}]=[0, T^{2}, \Lambda_{1,2}]$$
So $\Omega_{2}=T^{2}, T, 2T^{2}, 2T$ and we have $\widetilde{\Omega}_{2}=(1/2T) \Omega_{2}$, $a_{m(1,1)+1}=a_{2}=T$, $a_{m(2,1)+1}=a_{3}=2T^{2}$,
$a_{m(3,1)+1}=a_{4}=2T$ and we see that for all $1\leq k\leq 3$: $a_{4-k}=(2T)^{(-1)^{k+1}}a_{k +1}$.
$$R_{2,2}/S_{2,2}=[0, T^{2}, T, 2T^{2}, 2T, 2T^{2}, T^{4}]=[0, \Omega_{2,2}]=[0, \Omega_{2}, 2T^{2}, \Lambda_{2,2}]$$
\begin{eqnarray*}R_{3,2}/S_{3,2}&=&[0, T^{2}, T, 2T^{2}, 2T, 2T^{2}, T^{4}, 2T^{2}, 2T, 2T^{2}, T]=[0, \Omega_{3,2}]\\&=&
[0, \Omega_{2}, 2T^{2}, \Lambda_{2,2}, 2T^{2}, \Lambda_{3,2}]\end{eqnarray*}
$$R_{1,3}/S_{1,3}=[0, T^{2}, T, 2T^{2}, 2T, 2T^{2}, T^{4}, 2T^{2}, 2T, 2T^{2}, T, T^{2}]=[0, \Omega_{3}]=[0, T^{2}, \Lambda_{1,3}]$$
So $\Omega_{3}=\Omega_{1,3}=\Omega_{1,2}, 2T^{2}, \Lambda_{2,2}, 2T^{2}, \Lambda_{3,2}, T^{2}$ and we have
$a_{m(1,2)+1}=a_{5}=2T^{2}$, $a_{m(2,2)+1}=a_{7}=2T^{2}$, $a_{m(3,2)+1}=a_{11}=T^{2}$, $\Lambda_{2,2}=T^{4}=T^{-2}\Omega_{1}^{(3)}$,
$\Lambda_{3,2}=2T, 2T^{2}, T=\widetilde{\Lambda}_{1,2}$ and we see that for all $1\leq k\leq 10$: $a_{11-k}=a_{k +1}$.\\
\begin{eqnarray*}R_{2,3}/S_{2,3}&=&[0, T^{2}, T, 2T^{2}, 2T, 2T^{2}, T^{4}, 2T^{2}, 2T, 2T^{2}, T, T^{2}, T,
2T^{5}, 2T^{4}, T^{5}, T^{4}]\\&=&[0, \Omega_{2,3}]=[0, \Omega_{3}, T, \Lambda_{2,3}]\end{eqnarray*}
$$R_{3,3}/S_{3,3}=[0, T^{2}, T, 2T^{2}, 2T, 2T^{2}, T^{4}, 2T^{2}, 2T, 2T^{2}, T, T^{2}, T,
2T^{5}, 2T^{4}, T^{5}, T^{4},$$$$  2T^{2}, 2T, 2T^{2}, T, T^{2}, T, 2T^{5},
T, T^{2}, T, 2T^{2}]$$
$ \ \ \ \ \ \ \ \ \ \ \ \ \ \  \ \ \ \  \ \ =[0, \Omega_{3,3}]=[0, \Omega_{3}, 2T^{3}, \Lambda_{2,3},2T^{2}, \Lambda_{3,3}].$
$$R_{1,4}/S_{1,4}=[0, T^{2}, T, 2T^{2}, 2T, 2T^{2}, T^{4}, 2T^{2}, 2T, 2T^{2}, T, T^{2}, T,
2T^{5}, 2T^{4}, T^{5}, T^{4},$$$$ 2T^{2}, 2T, 2T^{2}, T, T^{2}, T, 2T^{5},
T, T^{2}, T, 2T^{2}, 2T]=[0, \Omega_{4}]$$
 So $\Omega_{4}=\Omega_{1,4}=\Omega_{1,3}, T, \Lambda_{2,3}, 2T^{2}, \Lambda_{3,3}, 2T$and we have
$a_{m(1,3)+1}=a_{12}=T$, $a_{m(2,3)+1}=a_{17}=2T^{2}$, $a_{m(3,3)+1}=a_{28}=2T$, $\Lambda_{2,3}=2T^{5}, 2T^{4}, T^{5}, T^{4}
=(1/2T)\Omega_{2}^{(3)}=2T\widetilde{\Lambda}_{2,3}$,
$\Lambda_{3,3}=2T, 2T^{2}, T, T^{2}, T, 2T^{5},
T, T^{2}, T, 2T^{2}=(1/2T)\widetilde{\Lambda}_{1,3}$ and we see that for all $1\leq k\leq 27$: $a_{28-k}=(2T)^{(-1)^{k+1}}a_{k +1}$.\\
\begin{rem}
Note that the equation $(W_{1})$ can be written as $\alpha=(C\alpha^{4}+1)/A$, so $\nu(\alpha^{4})=\nu(\alpha)$. Let $\beta=1/\alpha^{4}$. We will
determine the equation satisfied by $\beta$. We have $(A\alpha)^{4}=(C\alpha^{4}+1)^{4}=C^{4}\alpha^{16}+C^{3}\alpha^{12}+C\alpha^{4}+1$. Hence $\beta$
satisfies the equation $\beta^{4}+(-A^{4}+C)\beta^{3}+C^{3}\beta+C^{4}=0$. So it is clear that $\beta$ is hyperquadratic.
We can describe its continued fraction expansion as follow.
We put $\gamma=1/\alpha$. Then $\gamma$ satisfies the equation
\begin{eqnarray}\label{q}\gamma^{4}=A\gamma^{3}-C.\end{eqnarray}We know that the continued fraction expansion of $\gamma$ is
$$\gamma=[A, -A^{3}C^{-1},A^{3^{2}}C^{-2}, -A^{3^{3}}C^{-7},\ldots, (-1)^{n-1}A^{3^{n-1}}C^{-\frac{3^{n-1}+(-1)^{n}}{4}}, \ldots].$$
From the property \eqref{sd} of continued fractions and the equation \eqref{q} we get\begin{eqnarray*}
\gamma^{4}&=&[-C +A^{4}, -A^{3^{2}-1}C^{-3},A^{3^{3}+1}C^{-6}, -A^{3^{4}-1}C^{-21},\ldots, (-1)^{n-1}A^{3^{n}-(1)^{n}}C^{-\frac{3^{n}+3(-1)^{n}}{4}},
\ldots]\\&=&\beta=1/\alpha^{4}.
\end{eqnarray*}
This led us to deduce the following curious relation between square of continued fractions:
$$[0,A, -A^{3}C^{-1},A^{3^{2}}C^{-2}, -A^{3^{3}}C^{-7},\ldots, (-1)^{n-1}A^{3^{n-1}}C^{-\frac{3^{n-1}+(-1)^{n}}{4}}, \ldots]^{2}=[0, \Omega_{\infty}];$$
$$[0, \Omega_{\infty}]^{2}=[0,-C +A^{4}, -A^{3^{2}-1}C^{-3},A^{3^{3}+1}C^{-6},\ldots, (-1)^{n-1}A^{3^{n}-(1)^{n}}C^{-\frac{3^{n}+3(-1)^{n}}{4}},
\ldots].$$
\end{rem}
\begin{rem}
\noindent From Theorem \ref{t3}, the continued fraction expansion of $\beta$ solution of the equation $(W_{2})$ can be written as:
$$\beta=[0,A/C, (A/C)^{3}C,(A/C)^{3^{2}}C^{2}, \ldots, (A/C)^{3^{n-1}}C^{\frac{3^{n-1}+(-1)^{n}}{4}}, \ldots].$$
We wish to compute the continued fraction expansion and the approximation exponent of
$\alpha=\beta^{2}=[0,A/C, (A/C)^{3}C,(A/C)^{3^{2}}C^{2}, \ldots, (A/C)^{3^{n-1}}C^{\frac{3^{n-1}+(-1)^{n}}{4}}, \ldots]^{2}$. Note that,
from the equation $(W_{2})$, $\beta$ satisfies $\beta=(-\beta^{4}+C)/A$. So $\beta^{2}=(\beta^{4}-C)^{2}/A^{2}$, which gives that $\beta^{8}-2C\beta^{4}+C^{2}=A^{2}\beta^{2}$.
Then we deduce that $\alpha$ satisfies the equation
$$\alpha^{4}+C\alpha^{2}-A^{2}\alpha +C^{2}=0\ \ \ (E_{2}).$$
We have to state the following Conjecture.\\
\begin{conj}\label{tfg} Let $\alpha \in \mathbb{F}_{3}((T^{-1}))$ be the formal power series satisfying $(E_{2})$.\\
Let $(\Omega_{n})_{n\geq1}$ be a sequence of elements of $\mathbb{F}_{3}[T]$, defined by $\Omega_{0}=\emptyset$,\\
 $\Omega_{1}=A^{2}/C^{2}$ and for all $n\geq0$
\begin{eqnarray}\label{rrr}\left\{%
\begin{array}{ll}
    \Omega_{2n +1}= \Omega_{2n}, 2A^{2}/C^{2}, C^{2}\Omega_{2n-1}^{(3)}, 2A^{2}/C^{2}, \widetilde{\Omega}_{2n} \\
    \Omega_{2n +2}= \Omega_{2n +1}, 2A^{2}/C, C\Omega_{2n}^{(3)}, 2A^{2}/C^{2}, C\Omega_{2n +1} \\
\end{array}%
\right.    \end{eqnarray}
Then $\alpha=[0, \Omega_{\infty}]$.\end{conj}

Note that obtaining this conjecture was achieved in the same way as conjecture \ref{tf}. As it is so length, we omit it.
But it is interesting to state that we can add  the family of power series satisfying $(E_{2})$ to the set of elements admitting 2 as a value of
their approximation exponents agreeing with Roth value.

 \end{rem}
\noindent At the end, we point out a question related to this work:\\
\textbf{Open question:} Let $n$ be a positive integer and $(a_{i})_{1 \leq i \leq n}$ be a sequence of polynomials with coefficients in a finite field
such that $\deg a_i>0$. Let $D$ be a nonzero polynomial with coefficients in a finite field and with strictly positive degree such that $D$ divides
$a_{1}$. Suppose that $$[a_{n}, ....., a_{1}]=D^{-1}[a_{1}, ....., a_{n}].$$
Then $n$ is even and for all $0\leq k\leq n-1$:
$$a_{n-k}=\displaystyle\frac{1}{D^{(-1)^{k}}}a_{k +1}.$$

%\begin{acknowledgements}
%\end{acknowledgements}



\noindent Address of authors:
Department of Mathematics \\ % \hfill (Received 00 00 201?)\\
Sfax University   \\ %\hfill (Revised  00 00 201?)\\
Faculty of Sciences\\
Tunisia\\
Emails:{khalil.ayadi@isgis.usf.tn; azazaawatef91@gmail.com; beldisalah@gmail.com}

%% OTHER AUTHOR(S):
%\author[]{}
%\address{ }
%\email{}

%%%%%%%%%%%%%%  exercice &1    %%%%%%%%%%%%%%%%%%%%%%%%%%%%%%%%%%%%%%%%%%%%%%%%%%%%%%%%%%%%

% \ \newpage

% \setcounter{page}{1}
% \entete
%
% \Closesolutionfile{mycor}
% \Readsolutionfile{mycor}
%$\mathcal{U}(A)$ le groupe des éléments inversibles de $A$
\end{document}